\documentclass[11pt]{article}
\usepackage[margin=0.8in]{geometry}
\newcommand{\bfs}[1]{\mbox{\boldmath$#1$}}
\newcommand{\conv}[1]{\text{conv}\left({#1}\right)}
\newcommand{\ch}[1]{\text{ch}\left({#1}\right)}

\usepackage[ruled,linesnumbered,noline]{algorithm2e}
\usepackage[numbers,sort]{natbib}
\usepackage{amssymb}
\usepackage{mathtools}
\DeclarePairedDelimiter\ceil{\lceil}{\rceil}
\DeclarePairedDelimiter\floor{\lfloor}{\rfloor}
\usepackage{cases}
\usepackage{multirow}
\usepackage{algcompatible}
\usepackage{caption}
\providecommand{\keywords}[1]{\small \quad \quad \textbf{\textit{Keywords }} #1}
\usepackage{titlesec}
\titleformat{\section}
  {\normalfont\fontfamily{ptm}\fontsize{11}{11}\bfseries}{\thesection}{1em}{}
\titleformat{\subsection}
  {\normalfont\fontfamily{ptm}\fontsize{11}{11}\bfseries}{\thesubsection}{1em}{}
\titleformat{\subsubsection}
  {\normalfont\fontfamily{ptm}\fontsize{11}{11}\selectfont}{\thesubsubsection}{1em}{}
\newtheorem{theorem}{Theorem}[section]
\newtheorem{corollary}[theorem]{Corollary}
\newtheorem{lemma}[theorem]{Lemma}
\newtheorem{proposition}[theorem]{Proposition}
\newtheorem{definition}[theorem]{Definition}
\newtheorem{example}[theorem]{Example}
\usepackage{hyperref}
\hypersetup{colorlinks=true, urlcolor=cyan, linkcolor=cyan, citecolor=cyan}

\title{\large Mixed-Integer Programming Approaches to Generalized Submodular Optimization and its Applications}
\author{ \small Simge K\"u\c{c}\"ukyavuz \quad Qimeng Yu \vspace{0.2cm}  \\ \small Department of Industrial Engineering and Management Sciences \\ \small Northwestern University, Evanston, IL, USA \\ \small \{simge@northwestern.edu, kim.yu@u.northwestern.edu\}}
\date{\small \today} 

\begin{document}
\maketitle

\noindent {\small \textbf{Abstract.} Submodularity is an important concept in integer and combinatorial optimization. A classical submodular set function models the utility of selecting homogenous items from a single ground set, and such selections can be represented by binary variables. In practice, many problem contexts involve choosing heterogenous items from more than one ground set or selecting multiple copies of homogenous items, which call for extensions of submodularity. We refer to the optimization problems associated with such generalized notions of submodularity as Generalized Submodular Optimization (GSO). GSO is found in wide-ranging applications, including infrastructure design, healthcare, online marketing, and machine learning. Due to the often highly nonlinear (even non-convex and non-concave) objective function and the mixed-integer decision space, GSO is a broad subclass of challenging mixed-integer nonlinear programming problems. In this tutorial, we first provide an overview of classical submodularity. Then we introduce two subclasses of GSO, for which we present polyhedral theory for the mixed-integer set structures that arise from these problem classes. Our theoretical results lead to efficient and versatile exact solution methods that demonstrate their effectiveness in practical problems using real-world datasets.}\vspace{8pt} \\

\section{Introduction}
\label{sect:introduction}
Decision-making in complex systems commonly involves selecting items from a given collection to attain optimal utility. For example, in sensor placement, we are interested in selecting the sensor locations that provide maximal coverage. In such decision-making problems, we often observe the phenomenon of \emph{diminishing returns}. That is, the more items that we have selected, the lower is the marginal contribution of any additional item to the total utility. In the example of sensor placement, the total coverage of a sensor placement plan typically follows diminishing returns---as we install more sensors, the less additional coverage any newly added sensor contributes. The phenomenon of \emph{diminishing returns} is ubiquitous in real-world applications. Economies of scale, risk aversion, clustering, coverage, and influence propagation {all possess diminishing returns property}. 

\begin{figure}[h] 
   \centering
    \includegraphics[width=10cm]{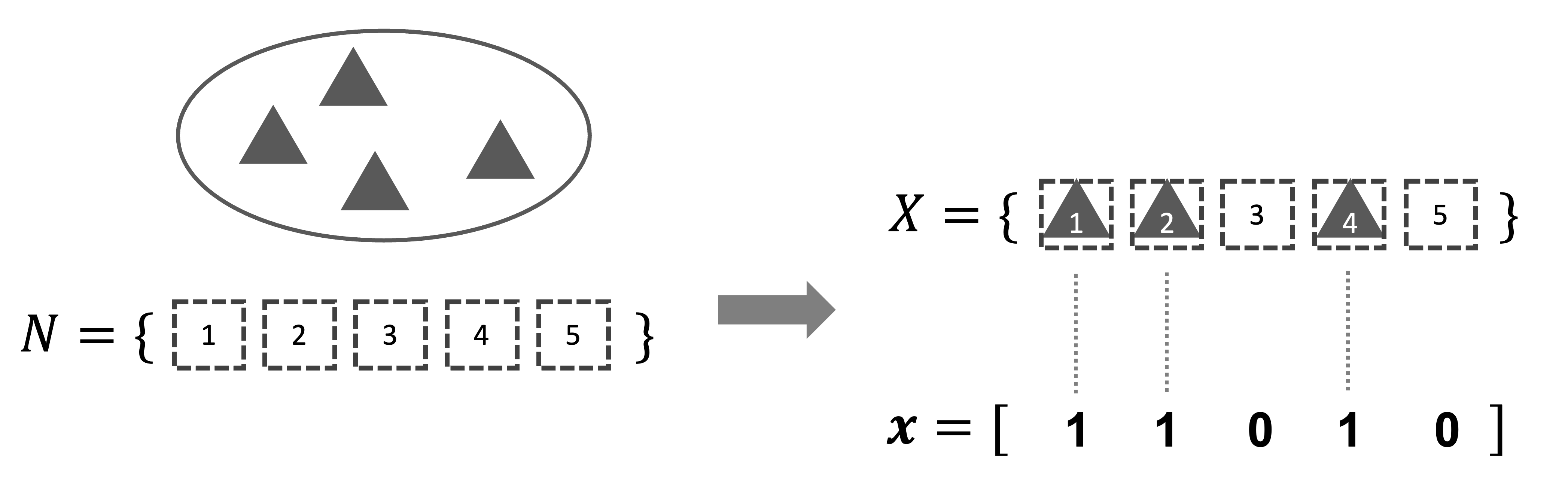}
    \caption{\small A pictorial illustration of submodular set function optimization. The triangles represent the given homogenous items, and $N$ is a finite ground set. The subset $X$ represents the decision of placing the items in the first, the second, and the fourth $\lq$bins'. By slightly abusing notation, the vector $\mathbf{x}$ provides a binary encoding of the same decision as $X$, with ones representing the locations at which we place an item, and zeros otherwise. }
    \label{fig:sub_set_opt}
\end{figure}

Mathematically, the concept of \emph{submodularity} formalizes diminishing returns. Submodularity is a property of \emph{set functions}, or functions defined over sets. Due to the prevalence of diminishing returns, submodular set functions are used to model utilities in a rich and diverse array of applications, such as facility location \cite{cornuejols1977uncapacitated}, image segmentation \cite{jegelka2011submodularity}, network influence propagation \cite{kempe2015maximizing,wu2018two}, and sensor placement \cite{krause2008near}. The associated decision-making problems are formulated as optimization problems with submodular set functions as their objectives, leading to \emph{submodular set function optimization} (submodular optimization, in short). In submodular optimization, we are given a supply of homogenous items (e.g., sensors in sensor placement) and a finite ground set (e.g., candidate sensor placement locations). We may think of the ground set as a set of $\lq$empty bins', and our goal is to determine a best subset of these $\lq$bins' in which we place the provided homogenous items. Whether we place an item in a $\lq$bin' or not is naturally modeled by a binary decision variable, making submodular set function optimization a class of discrete optimization problems. Figure \ref{fig:sub_set_opt} provides a pictorial high-level illustration of submodular optimization. 

\begin{figure}[h] 
   \centering
    \includegraphics[width=10cm]{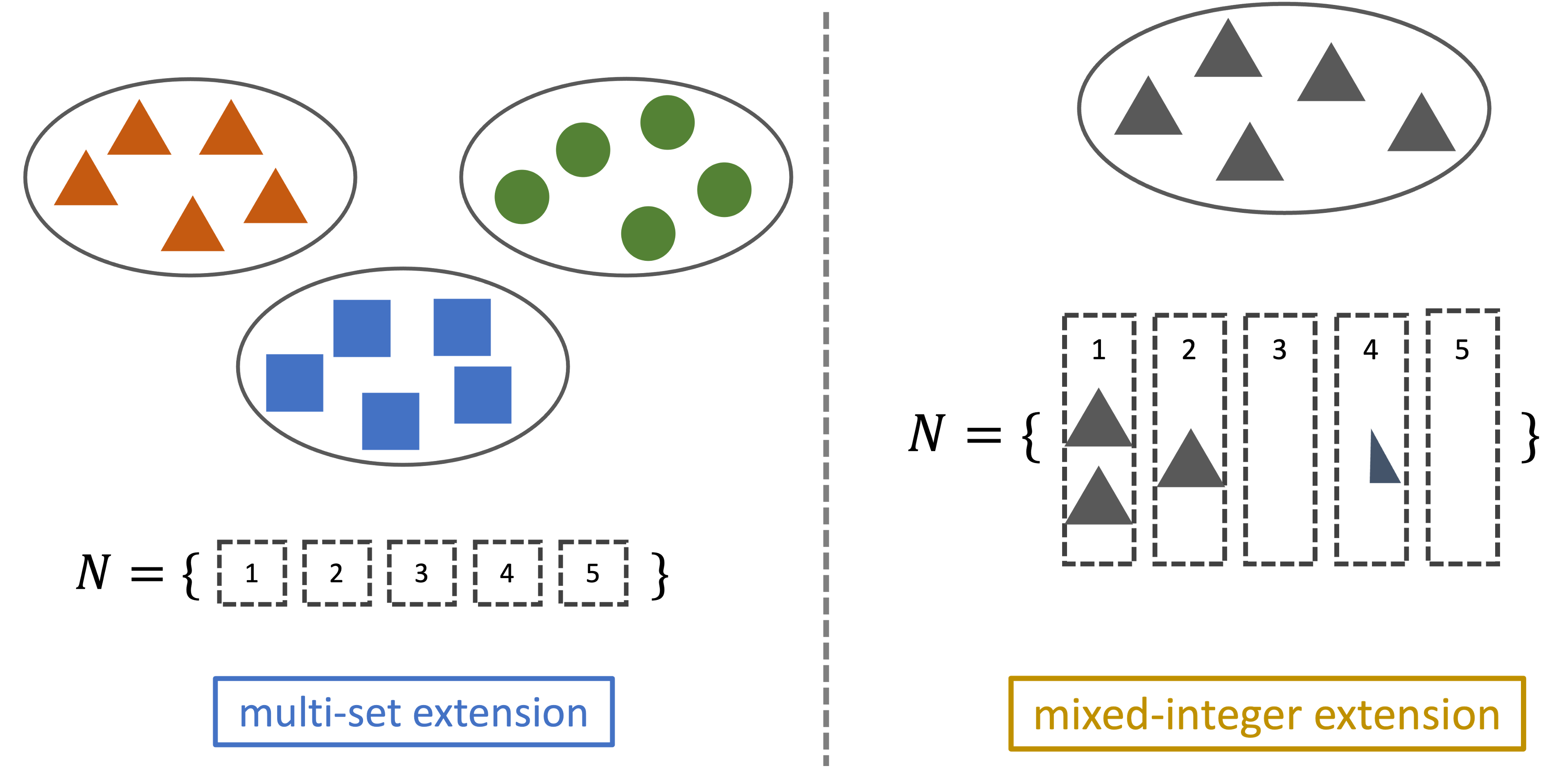}
    \caption{\small Two scenarios that call for generalizations of classical submodularity. The first scenario on the left calls for the multi-set extension. The second scenario on the right gives rise to the mixed-integer extension.}
    \label{fig:sub_extensions}
\end{figure}

Despite the wide-ranging modeling power of submodular set functions, applications increasingly call for more general modeling tools that still capture diminishing returns. Consider the two scenarios depicted in Figure \ref{fig:sub_extensions}. In the first scenario on the left, we are given a supply of \emph{heterogenous} items instead. Now our task is not only to determine the best subset of the $\lq$bins' to fill, but also the best type of items for the selected $\lq$bins'. This calls for the multi-set extension of submodularity. 
In the second scenario on the right, we are given a finite ground set of $\lq$bins' and a homogenous collection of items. However, we are now allowed to place multiple, and even fractional, copies of the items in each $\lq$bin'. This scenario gives rise to the mixed-integer extension of submodularity. We refer to classical submodularity and all its extensions by the unifying term \emph{Generalized Submodularity} (GS).
We summarize two examples of GS below. More applications are described in detail in Sections \ref{sect:multiset} and \ref{sect:mi}. 

\textbf{Multi-type sensor placement.} In practice, we may be provided with multiple types of sensors that monitor different aspects of the environment (e.g., temperature, humidity, and luminance), and each candidate sensor location can hold one sensor. We now need to decide at which locations we place a sensor, and additionally, which type of sensor we install at each selected location. The multi-set extension of submodularity is an ideal modeling tool for this multi-type sensor placement problem. 

\textbf{Sensor energy management.} In another sensor placement application, we are given sensors that monitor the same environmental factors but with adjustable energy levels. The energy levels could be continuous ranges or discrete dials, and our goal is to assign the energy levels that yield optimal utility. This problem fits under the second scenario that calls for the mixed-integer extension. 

In this tutorial, we introduce notions of GS and demonstrate their modeling power with real-world applications. In addition, we review mixed-integer programming-based (polyhedral) approaches to optimize GS functions, with or without constraints, to global optimality. 
The polyhedral approach has been highly successful in attaining global optimal solutions to classical submodular optimization under complicating constraints \cite{edmonds2003submodular,wolsey1999integer,ahmed2011maximizing,yu2017maximizing,shi2022sequence,yu2017polyhedral,yu2023strong, Coniglio2022} as well as for general set functions \citep{atamturk2021submodular}. Furthermore, submodularity has been uncovered and exploited in improving the formulations of mixed-binary convex quadratic and conic optimization problems \cite{gomez2018submodularity, atamturk2020submodularity, atamturk2020supermodularity, kilincc2020conic, atamturk2008polymatroids, atamturk2019lifted}. Furthermore, these approaches have been fruitful in stochastic and risk-averse problems that involve submodular structures \citep{wu2018two, wu2019probabilistic, wu2020exact, kilincc2021joint, xie2019distributionally, zhang2018ambiguous, shen2022chance}. We focus our attention to expanding the reach of these approaches to generalized submodular functions.

The rest of this tutorial is structured as follows. We first provide preliminaries on classical submodularity in Section \ref{sect:preliminaries}. We demonstrate how classical submodularity can be exploited to facilitate solution of a class of constrained submodular set function minimization problems in Section \ref{sect:concave_min}. In Sections \ref{sect:multiset} and \ref{sect:mi}, we formally state two generalizations of submodularity and summarize their applications. We present the polyhedral results for the key mixed-integer sets arising from the associated GSO, with which we describe exact solution methods for these optimization problems. We conclude this tutorial with Section \ref{sect:conclusion}.

\section{Preliminaries} 
\label{sect:preliminaries}
In this section, we provide an overview of classical submodularity. Let $N = \{1,2,\dots,n\}$ denote a finite non-empty ground set, and let $2^N$ be the collection of all the subsets of $N$. Formally, submodular set functions are defined as follows. 

\begin{definition}
\label{def:submodular}
A function $f:2^N \rightarrow \mathbb{R}$ is \emph{submodular} if 
\begin{equation}
\label{eq:sub_def1}
f(X) + f(Y) \geq f(X\cap Y) + f(X\cup Y)
\end{equation}
for any $X, Y\subseteq N$. 
\end{definition}

A function $f$ is \emph{supermodular} if inequality \eqref{eq:sub_def1} is reversed for all $X, Y\subseteq N$, or equivalently, if $-f$ is submodular. Moreover, function $f$ is \emph{monotone non-decreasing}, or \emph{monotone} for brevity, if $f(X) \leq f(Y)$ for all $X\subseteq Y$. For any subset $X\subset N$ and any item $i\in N$, we define the \emph{marginal return} of adding $i$ to $X$ as 
\begin{equation}
\label{eq:rho}
\rho_i(X) := f(X\cup\{i\}) - f(X).
\end{equation}
Intuitively, submodularity is equivalent to diminishing returns. The next alternative definition for submodular set functions clearly reflects this intuition. 
\begin{figure}[h] 
   \centering
   \includegraphics[width=9cm]{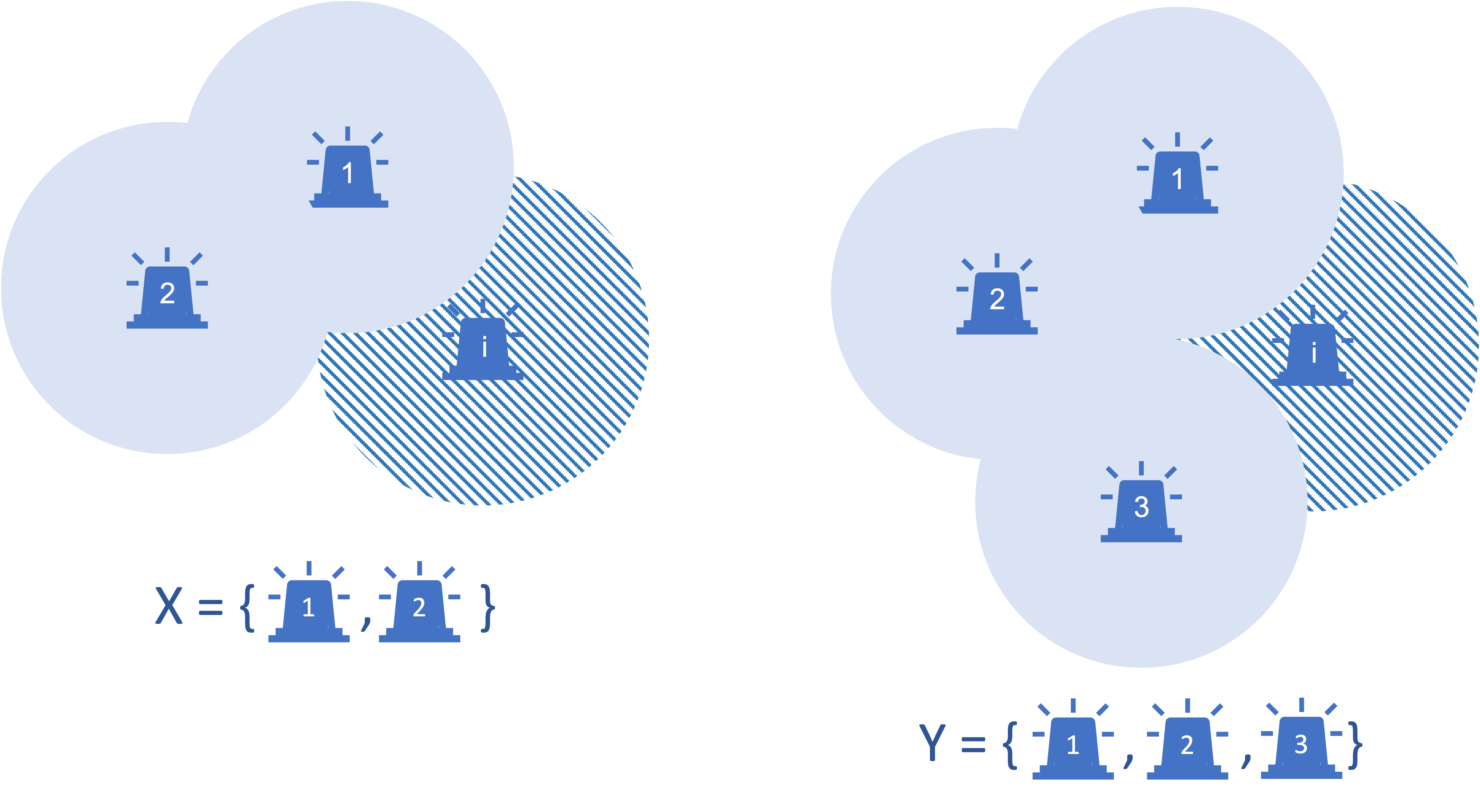}
   \caption{\small An illustration of diminishing returns in sensor placement. Each circle represents the coverage of each sensor. The sets $X$ and $Y$ are two sets of locations at which we have installed sensors. The solid regions correspond to the total coverages of placement plans $X$ and $Y$, respectively. The striped regions highlight the marginal contributions of an additional sensor installed at location $i$ to the two plans. The more sensors that we have installed, the less additional coverage we will obtain by installing more sensors. }
   \label{fig:sensor_sub}
\end{figure}  
\begin{definition}
\label{def:submodular_equiv}
A function $f:2^N \rightarrow \mathbb{R}$ is \emph{submodular} if 
\begin{equation}
\rho_i(X) \geq \rho_i(Y)
\end{equation}
for any $X\subseteq Y\subseteq N$ and any $i\in N\backslash Y$. 
\end{definition}
By Definition \ref{def:submodular_equiv}, the marginal return of item $i$ diminishes as we switch from selection $X$ to an inclusion-wise larger selection, $Y$. Figure \ref{fig:sensor_sub} demonstrates diminishing returns in sensor placement. We may determine whether a set function is submodular and monotone non-decreasing using the following lemma. 
\begin{lemma} \cite{wolsey1999integer}
A function $f:2^N \rightarrow \mathbb{R}$ is \emph{submodular and monotone} if and only if for all $X, Y\subset N$, 
\[f(Y)\leq f(X) +\sum_{i\in Y\setminus X} \rho_i(X).\]
\end{lemma}

\textbf{Examples of submodular set functions.} Next we provide a few examples of submodular set functions which frequently appear in combinatorial optimization. Recall $N = \{1,2,\dots,n\}$ denotes a finite non-empty ground set. We consider the following set functions $f:2^N\rightarrow \mathbb{R}$. 
\begin{itemize}
\item The \emph{linear} function $f(X)=\sum_{i\in X} a_i$, where $\mathbf{a}\in\mathbb{R}^n$, is both submodular and supermodular. Thus, we say $f$ is \emph{modular}. If $a_i\geq 0$ for all $i\in N$, then $f$ is also monotone. 
\item Suppose we are given a set of items, $\mathcal{C}$, and $n$ subsets of these items, $A_1, \dots, A_n \subset \mathcal{C}$. The \emph{set coverage} function $f(X) = |\bigcup_{i\in X} A_i|$ is submodular.   
\item Let $\mathcal{G} = (N, \mathcal{A})$ be a directed graph with vertex set $N$ and arc set $\mathcal{A}$. In addition, every arc $(i,j)\in\mathcal{A}$ has a capacity $c_{ij}\geq 0$. The \emph{graph cut} function $f(X)=\sum_{(i,j)\in \mathcal{A}: i\in X, j\in  N\setminus X} c_{ij}$ is submodular but not monotone. 
\end{itemize}

\textbf{Operations on submodular set functions.} The following list includes some operations that preserve submodularity. Note that the functions $h:2^N\rightarrow \mathbb{R}$ below are set functions defined over $N = \{1,2,\dots,n\}$. 
\begin{itemize}
\item Let $f:2^N\rightarrow \mathbb{R}$ be a submodular set function. For any constant $c\in\mathbb{R}$, the function $h(X) := f(X) - c$, for all $X\subseteq N$, is submodular. 
\item Suppose $f, g:2^N\rightarrow \mathbb{R}$ are two submodular set functions. Then for any nonnegative constants $\alpha, \beta\in\mathbb{R}_+$, the function $h(X) := \alpha f(X) +\beta g(X)$, for all $X\subseteq N$, is submodular. More generally, any conic combination of submodular functions is submodular. 
\item Let $f:\mathbb{R}\rightarrow\mathbb{R}$ be a concave function and $g:2^N\rightarrow \mathbb{R}$ be a nonnegative linear function. The composition $h(X) := f(g(X))$, for all $X\subseteq N$, is submodular. 
\item Suppose $f, g:2^N\rightarrow \mathbb{R}$ are two submodular set functions, such that $f-g$ is monotone non-decreasing or non-increasing. The minimum $h(X) := \min\{f(X), g(X)\}$ is submodular. 
\end{itemize}

We provide the definitions and properties above to help the readers identify submodularity in applications. To illustrate how submodularity can be exploited to attain optimal decisions in these applications, we further summarize a few key results of submodular set function optimization. 

Throughout the remaining discussion, we assume that $f(\emptyset) = 0$. This assumption is without loss of generality because otherwise we can equivalently consider $f'(X) := f(X) - f(\emptyset)$ for all $X\subseteq N$. If $f$ is submodular and/or monotone, then so is $f'$. We also assume that $f$ is available explicitly via a value oracle---given any $X\subseteq N$, $f(X)$ is returned in time that is polynomial in $n$. By slightly abusing notation, we use $X\subseteq N$ and its binary representation interchangeably. To be precise, any $X\subseteq N$ uniquely corresponds to a binary vector $\mathbf{x}\in\{0,1\}^n$, where 
\[x_i = \begin{cases}
1, & \text{if } i\in X, \\
0, & \text{otherwise.}
\end{cases}\]
Now a set function $f:2^N\rightarrow\mathbb{R}$ can be rewritten as a function defined over binary vectors, $f: \{0,1\}^n \rightarrow \mathbb{R}$.

\subsection{Submodular set function minimization}
We first consider the \emph{unconstrained} minimization problem  
\begin{equation}
\label{eq:uncon_sub_min}
\min_{X\subseteq N} f(X),
\end{equation}
where $f$ is a submodular set function. This class of optimization problems is known to be strongly polynomial-time solvable \cite{grotschel1981ellipsoid, iwata2001combinatorial, lee2015faster, orlin2009faster, cunningham1985submodular, schrijver2000combinatorial, mccormick2005submodular, iwata2009simple, iwata2008submodular}. In fact, $f$ can be further extended to the continuous domain $[0,1]^n$ via \emph{Lov\'asz extension} \cite{lovasz1983submodular}. 

\begin{definition}
\label{def:lovasz}
Suppose we are given an arbitrary function $f:\{0,1\}^n\rightarrow \mathbb{R}$. The \emph{Lov\'asz extension} $f^L:[0,1]^n\rightarrow \mathbb{R}$ of $f$ is defined by
 \begin{equation}
 \label{eq:lovasz}
 f^L(\mathbf{x})=\sum_{i=1}^n (x_{\delta(i)}-x_{\delta(i+1)})f\left(X^{\bfs{\delta},i}\right),
 \end{equation}
for any $\mathbf{x}\in [0,1]^n$. Here, $\bfs{\delta} = (\delta(1), \delta(2), \dots, \delta(n))$ is a permutation of $N = \{1,2, \dots, n\}$, such that $x_{\delta(1)} \geq x_{\delta(2)} \geq \dots \geq x_{\delta(n)}$. We let $x_{\delta(n+1)}=0$, and $X^{\bfs{\delta}, i} := \{\delta(1), \dots, \delta(i)\}$ for any $i\in N$. 
\end{definition}

We note that \eqref{eq:lovasz} is equivalent to
\[f^L(\mathbf{x})=\sum_{i=1}^n \left[f\left(X^{\bfs{\delta},i}\right)-f\left(X^{\bfs{\delta},i-1}\right)\right] x_{\delta(i)} = \sum_{i=1}^n \rho_{\delta(i)}\left(X^{\bfs{\delta},i-1}\right) x_{\delta(i)},\]
where $X^{\bfs{\delta}, 0}=\emptyset$, and $\rho_\cdot(\cdot)$ is defined as in \eqref{eq:rho}. While $f^L$ appears to be linear, due to the fact that the permutation $\bfs{\delta}$ depends on $\mathbf{x}$, this is not necessarily the case. Nevertheless, the Lov\'asz extension is easy to compute given an efficient value oracle. When $f$ is submodular, $f^L$ has the following desirable property. 
\begin{proposition} \cite{lovasz1983submodular}
The Lov\'asz extension of $f$ coincides with its convex closure $f^c$  if and only if $f$ is a submodular set function, where 
\[ f^c(\mathbf{x})=\min\left\{\sum_{S\subseteq N} \alpha_S f(S): \sum_{S\subseteq N} \alpha_S \mathbf 1_S=\mathbf{x},  \sum_{S\subseteq N} \alpha_S=1, \alpha_S\ge 0, \forall \; S\subseteq N\right\},\]  and $\mathbf{1}_S$ is a binary vector with ones at $i\in S$ and zeros in all other entries. 
\end{proposition}
In other words, minimizing the convex closure of the submodular function $f$ yields binary optimal solutions. Therefore, unconstrained submodular set function minimization is equivalent to minimizing a convex function, which suggests that this problem can be solved efficiently. 

The complexity of unconstrained submodular set function minimization also follows from the polynomial equivalence of optimization and separation. We next discuss this class of minimization problems from a polyhedral perspective. Let 
\[\mathcal{Q}_f := \{(\mathbf{x}, w)\in\{0,1\}^n \times \mathbb{R} : w \geq f(\mathbf{x})\}\]
be the epigraph of $f$. Problem \eqref{eq:uncon_sub_min} has an equivalent reformulation: 
\begin{equation}
\label{eq:uncons_sub_min_reform}
\min \{w : (\mathbf{x}, w)\in\conv{\mathcal{Q}_f}\}.
\end{equation}
The convex hull of the epigraph, $\conv{\mathcal{Q}_f}$, is described by the trivial inequalities $\bfs{0}\leq \mathbf{x}\leq \bfs{1}$ and all the \emph{extended polymatroid inequalities} (EPI) \cite{lovasz1983submodular, atamturk2008polymatroids, atamturk2021submodular}. Each EPI is associated with a permutation of $N$ and is defined as follows. 
\begin{definition}
Given any permutation $\bfs{\delta}$ of $N$, the corresponding \emph{extended polymatroid inequality (EPI)} is 
\begin{equation}
\label{eq:EPI}
w \geq \sum_{i=1}^n \rho_{\delta(i)}\left(X^{\bfs{\delta},i-1}\right) x_{\delta(i)}, 
\end{equation}
where $X^{\bfs{\delta},\cdot}= \{\delta(1), \dots, \delta(i)\}$ for any $i\in N$, $X^{\bfs{\delta}, 0}=\emptyset$, and $\rho_\cdot(\cdot)$ is defined as in \eqref{eq:rho}. 
\end{definition}
The convex hull $\conv{\mathcal{Q}_f}$ in formulation \eqref{eq:uncons_sub_min_reform} is defined by all the EPIs and $\bfs{0}\leq \mathbf{x}\leq \bfs{1}$, making this formulation a linear program with exponentially many constraints. 

\textbf{Delayed constraint generation framework.} To tackle optimization problems with large numbers of constraints, it is natural to apply the \emph{Delayed Constraint Generation (DCG)} framework. Consider reformulation \eqref{eq:uncons_sub_min_reform} as an example. Instead of entering this full master problem with exponentially many EPIs into a solver, we can iteratively solve \emph{relaxed master problems} that have fewer constraints and add necessary constraints in each iteration. A relaxed master problem assumes the form of 
\begin{equation}
\label{eq:dcg_relaxed}
\begin{aligned}
\min & \quad w  \\
\text{s.t.} & \quad (\mathbf{x},w)\in \mathcal{C}, 
\end{aligned}
\end{equation} 
where $\mathcal{C}$ is constructed by only a subset of the constraints defining $\conv{\mathcal{Q}_f}$ in Problem \eqref{eq:uncons_sub_min_reform}. Algorithm \ref{alg:general_dcg} provides an outline of the DCG algorithm. The term $(\text{UB}-\text{LB})/\text{UB}$ measures the optimality gap, where LB and UB denote the lower and upper bounds on the objective value, respectively. In this algorithm, while the optimality gap exceeds a pre-specified tolerance parameter $\epsilon$, we solve the relaxed master problem \eqref{eq:dcg_relaxed} to obtain its optimal solution and objective pair $(\overline{\mathbf{x}}, \overline{w})$. Given that \eqref{eq:dcg_relaxed} is a relaxed version of the original problem \eqref{eq:uncons_sub_min_reform}, $\overline{w}$ serves as a lower bound on the optimal objective. The function evaluation $f(\overline{\mathbf{x}})$ is an upper bound because $\overline{\mathbf{x}}$ is feasible. If $\overline{w}$ does not match $f(\overline{\mathbf{x}})$, then we add a violated valid inequality to $\mathcal{C}$ (line 7). The problem of determining such an inequality is called a \emph{separation problem}. Ideally, we would like to find a \emph{most violated valid inequality}, which is the most powerful constraint that will cut off $(\overline{\mathbf{x}}, \overline{w})$ in the next iteration of DCG. {Furthermore, by aiming to maximize violation, we are also able to identify if no violation exists.} We solve the updated relaxed master problem and repeat until the optimality gap is sufficiently small. 

\vspace{0.2cm}
\begin{algorithm}[H]
 \caption{\texttt{Delayed\_Constraint\_Generation}}
 \label{alg:general_dcg}
\begin{algorithmic}[1]
\STATE \textbf{Input} initial $\mathcal{C}$, $\text{LB} = -\infty$, $\text{UB} = \infty$\;
\WHILE  {$(\text{UB}-\text{LB})/\text{UB}>\epsilon$}
    \STATE Solve the current relaxed master problem \eqref{eq:dcg_relaxed} to get $(\overline{\mathbf{x}}, \overline{w})$\;
    \STATE $\text{LB}\leftarrow\overline{w}$\;
    \STATE Compute $f(\overline{\mathbf{x}})$\;
    \IF {$\overline{w}<f(\overline{\mathbf{x}})$}
        \STATE Determine a violated valid inequality $w\geq \pi^\top x$ and add it to update $\mathcal{C}$\;
    \ENDIF
     \IF {$\text{UB}>f(\overline{\mathbf{x}})$}
        \STATE $\text{UB}\leftarrow f(\overline{\mathbf{x}})$\;
    \ENDIF
    \STATE Update the incumbent solution to $\overline{\mathbf{x}}$\;
\ENDWHILE
\STATE  \textbf{Output} $\overline{w}$, $\overline{\mathbf{x}}$. 
\end{algorithmic}
\end{algorithm}

\vspace{0.4cm} As indicated by the Ellipsoid method, informally, optimization problems $\min\{\bfs{c}^\top\mathbf{x} : \mathbf{x}\in\mathcal{X}\}$ are polynomially solvable if the separation problem for $\conv{\mathcal{X}}$ can be solved in polynomial time. The separation problem for $\conv{\mathcal{Q}_f}$ in Problem \eqref{eq:uncons_sub_min_reform} happens to be $\lq$easy'. The most violated EPI for any $\mathbf{x}\in[0,1]^n$ can be found using an $\mathcal{O}(n\log n)$ greedy algorithm proposed by \cite{edmonds2003submodular}. 

In contrast to the unconstrained case, constrained submodular minimization (e.g., with cardinality constraints) is generally NP-hard \cite{svitkina2011submodular}. There exist exceptions when the problems satisfy special properties, such as minimizing the composition of a concave function with a non-negative affine function under a cardinality constraint \cite{hassin1989maximizing, onn2003convex}. We discuss the mixed-integer programming approach to this particular problem class in Section \ref{sect:concave_min}. To attain the exact optimal solutions for general constrained submodular minimization, we may implement reformulation \eqref{eq:uncons_sub_min_reform} with additional linearly representable constraints in the DCG framework (see \cite{atamturk2008polymatroids} for example).

\subsection{Submodular set function maximization}
Submodular maximization is NP-hard even in the unconstrained case. Given the complexity, many heuristics and approximation algorithms are proposed \cite{nemhauser1978best, lee2010maximizing, calinescu2007maximizing, sviridenko2004note, orlin2018robust}. In particular, \cite{nemhauser1978analysis} propose a celebrated  $(1-1/e)$-approximation algorithm for the cardinality-constrained maximization problem:
\[ \max_{X:|X|\leq k} f(X),\]
when $f$ is monotone and submodular. Here, $e$ is the base of normal logarithm. This algorithm is a greedy method that sequentially collects $k$ items with the highest marginal returns. 

In many applications, especially those that involve high-capital investments like facility location and sensor network design, it is desirable to solve for the exact optimal solutions rather than resorting to suboptimal ones. This is because the exact optimal solutions can lead to significant utility improvements and cost reductions. In \cite{wolsey1999integer}, the authors tackle submodular maximization, namely 
\begin{equation}
\label{eq:sub_max}
\max_{X\subseteq N} f(X), 
\end{equation} 
using an integer programming approach. We note that the submodular objective function $f$ does not have to be monotone. Let $f^*(X) := f(X) - \sum_{j\in X} \rho_j(N\setminus\{j\})$ for all $X\subseteq N$. The authors propose the following equivalent mixed-integer programming formulation:
\begin{equation}
\label{eq:sub_max_mip}
\max_{ (\mathbf{x}, w) \in \{0,1\}^n\times \mathbb{R} } \{w : (\mathbf{x}, w)\in \mathcal{Q}_f\},
\end{equation}
where 
\[\mathcal{Q}_f := \{(\mathbf{x}, w)\in \mathbb{R}_+^n\times \mathbb{R}: w \leq f^*(S)+\sum_{j\in N\setminus S} [f^*(S\cup \{j\})-f^*(S)] x_j + \sum_{j\in N} \rho_j(N\setminus \{j\}) x_j, \forall \; S\subseteq N\}\]
{is the hypograph of $f$ when $\mathbf{x}$'s are restricted to binary vectors.} 
Each inequality 
\begin{equation}
\label{eq:sub_ineq}
w \leq f^*(S)+\sum_{j\in N\setminus S} [f^*(S\cup \{j\})-f^*(S)] x_j + \sum_{j\in N} \rho_j(N\setminus \{j\}) x_j
\end{equation}
is called a \emph{submodular inequality} associated with $S\subseteq N$. Formulation \eqref{eq:sub_max_mip} has exponentially many constraints, which can be solved by DCG. The mixed-integer programming approach is able to handle additional constraints, and Formulation \eqref{eq:sub_max_mip} has been further strengthened for constrained submodular maximization by \cite{ahmed2011maximizing,yu2017maximizing,shi2022sequence}.

\section{Cardinality-Constrained Submodular Minimization}
\label{sect:concave_min}
To illustrate how to exploit submodularity to facilitate decision-making and how to apply the mixed-integer programming approach, this section focuses on a class of constrained submodular set function minimization problems. 

As mentioned in the previous section, the composition of a univariate concave function and a nonnegative linear function is a submodular set function. That is, for any $\mathbf{a}\in\mathbb{R}^n_+$ and any concave function $f:\mathbb{R}\rightarrow\mathbb{R}$, $F(X) := f\left(\sum_{i\in X} a_i \right)$ for all $X\subseteq N$, is submodular. Alternatively, $F$ can be written as $F(\mathbf{x}) = f\left(\sum_{i=1}^n a_i x_i\right)$ for all $\mathbf{x}\in\{0,1\}^n$. 

Such functions are frequently used to model utilities in problems that involve risk aversion or economies of scale, including mean-risk optimization  \cite{atamturk2008polymatroids, atamturk2019lifted} and concave cost facility location \cite{feldman1966warehouse, hajiaghayi2003facility}. Mean-risk optimization is a powerful modeling tool in financial applications, including investment portfolio management. We again observe diminishing returns in these applications. In investment portfolio management, for example, an investor's perception of risk can be modeled using the aforementioned functions due to the diminishing returns intuition---as an investor's portfolio becomes more aggressive, she tends to perceive the risk associated with any additional financial product to be milder. The number of financial products chosen is limited by a cardinality upper bound, $k$. Therefore, it is of practical interest to consider the following problem: 
\begin{equation}
\label{prob:cc_concave_sub_min}
\min \left\{  f\left(\sum_{i=1}^n a_ix_i\right) : \mathbf{x}\in\{0,1\}^n, \sum_{i=1}^n x_i \leq k \right\}.
\end{equation}
 Problem \eqref{prob:cc_concave_sub_min} has the following reformulation: 
\[\min \left\{ w : (\mathbf{x}, w)\in\mathcal{P}_k^m\right\},  \]
where
\begin{equation}
\label{eq:P}
\mathcal{P}_k^m = \left\{ (\mathbf{x}, w)\in\times\{0,1\}^n \times \mathbb{R} : w \geq f\left(\sum_{i=1}^n a_ix_i\right), \sum_{i=1}^n x_i \leq k \right\}
\end{equation}
represents the epigraph of $f$ under the cardinality constraint. The superscript $m\in\{1,2,\dots, n\}$ denotes the number of distinct values in $\mathbf{a}$. Problem \eqref{prob:cc_concave_sub_min} is polynomially solvable \cite{hassin1989maximizing, onn2003convex}, despite the NP-hardness of general constrained submodular minimization. This complexity result suggests that a complete convex hull characterization of $\mathcal{P}_k^m$ may be tractable. Such a characterization is desirable because it leads to a versatile branch-and-cut approach to handle more constraints in addition to the cardinality constraint. However, a tractable characterization for $\conv{\mathcal{P}_k^m}$ has been an open problem. 

Recall from Section \ref{sect:preliminaries} that the EPIs \eqref{eq:EPI} describe the epigraph convex hull without the cardinality constraint. The authors in \cite{yu2017polyhedral} account for the cardinality constraint and provide the full description of $\conv{\mathcal{P}_k^1}$ (i.e., when $a_i = \alpha \geq 0$ for all $i\in N$). The authors propose strong valid linear inequalities called the \emph{separation inequalities (SIs)}. 
\begin{definition}
\label{def:SI}
For any permutation $\bfs{\delta} = (\delta(1),\delta(2),\dots, \delta(n))$ of $N$ and a fixed integer parameter $0\leq i_0\leq  k-1$, a \emph{separation inequality (SI)} is given by  
\begin{equation}
\label{eq:sepa_coeff}
w\geq \sum_{i=1}^{i_0} \rho_{\delta(i)}x_{\delta(i)} + \sum_{i=i_0+1}^{n} \psi x_{\delta(i)}, 
\end{equation} 
where $\psi := \frac{f(k \alpha) - f(i_0 \alpha)}{k-i_0}$, and $\rho_{\delta(i)}:= f(i \alpha) - f((i-1)\alpha)$. 
\end{definition}
The authors show that the SIs, the cardinality constraint, and $\bfs{0}\leq \mathbf{x}\leq \bfs{1}$, completely describe $\conv{\mathcal{P}^1_k}$. Furthermore, they obtain a class of valid inequalities for $\conv{\mathcal{P}^m_k}$ with general $m$, given by Definition \ref{def:ALI}. For ease of notation, given any permutation $\bfs{\delta}$ of $N$, we re-index $N$ such that $\bfs{\delta}$ is the natural order $(1,2,\dots,n)$. Thus we omit $\bfs{\delta}$ in the notation below. 
\begin{definition}
\label{def:ALI}
The \emph{approximate lifted inequalities}, or ALIs, assume the form of
\[ w\geq \sum_{i=1}^k\rho_i x_i + \sum_{i=k+1}^n \phi_ix_i. \]
The coefficients $\rho_i$ for $i\in\{1,\dots, k\}$ are consistent with those in \eqref{eq:sepa_coeff}. For each $i>k$, let $T$ with $|T|= k-1$ be a subset of $\{1,\dots,i-1\}$ such that the sum of the weights are as high as possible, and $\phi_i = f(a_i + \sum_{j\in T}a_j) - f(\sum_{j\in T}a_j)$.
\end{definition}
Stronger valid inequalities have also been proposed for the general case (i.e., arbitrary $m$) without explicit forms of the coefficients.

In \cite{yu2023strong}, we make further progress in the case where $\mathbf{a}$ contains up to two distinct values and show multiple ways to apply their results to the general case (i.e., any $m\geq 1$). When $m=2$, we denote the two weights in $\mathbf{a}$ by $a_L$ and $a_H$, respectively, with $a_L \leq a_H$. We also define $\mathcal{I}_L = \{i\in N : a_i = a_L\}$ and $\mathcal{I}_H = \{i\in N : a_i = a_H\}$ as the index sets corresponding to the lower and higher values in $\mathbf{a}$. The results that we describe next hold for any permutation $\bfs{\delta}$ of $N$. The ground set $N$ is re-indexed so that $\bfs{\delta}$ is the natural order and thus is dropped from the notation. Consider any subset $S\subseteq N$ with cardinality $k$. By re-indexing, we assume $S = \{1,\dots, k\}$ without loss of generality. 

We now share two important observations. First, consider the set
\[\mathcal{P}^2_k(S) = \left\{ (\mathbf{x}, w)\in\{0,1\}^k \times \mathbb{R}: w \geq f\left(\sum_{i=1}^k a_ix_i\right) \right\}, \]
which is exactly $\mathcal{P}^2_k$ with $x_i$ fixed to zeros for all $i\in N\setminus S$. The cardinality constraint trivially holds because at most $k$ variables could be non-zero. We observe that each EPI \eqref{eq:EPI}, $w \geq \sum_{i=1}^k \rho_i x_i$, associated with any ordering of $S$ is facet-defining (i.e., among the strongest valid inequalities) for $\mathcal{P}^2_k(S)$. Next, consider the sets $\mathcal{P}^1_k(\mathcal{I}_L)$ and $\mathcal{P}^1_k(\mathcal{I}_H)$, where
\begin{equation}
\label{eq:L_P1k}
\mathcal{P}^1_k(\mathcal{I}_L) = \left\{ (w, x)\in\mathbb{R}\times\{0,1\}^{\mathcal{I}_L} : w \geq f\left(a_L \sum_{i\in\mathcal{I}_L} x_i\right), \sum_{i\in\mathcal{I}_L} x_i \leq k \right\},
\end{equation} and 
\begin{equation}
\label{eq:H_P1k}
\mathcal{P}^1_k(\mathcal{I}_H) = \left\{ (w, x)\in\mathbb{R}\times\{0,1\}^{\mathcal{I}_H} : w \geq f\left(a_H \sum_{i\in\mathcal{I}_H} x_i\right), \sum_{i\in\mathcal{I}_H} x_i \leq k \right\}. 
\end{equation}
In $\mathcal{P}^1_k(\mathcal{I}_L)$, all the variables $x_i$ with $i\in \mathcal{I}_H$ are fixed to zeros, and $\mathcal{P}^1_k(\mathcal{I}_H)$ imposes the same restrictions on $i\in \mathcal{I}_L$. We note that the SIs are facet-defining for these two sets. These two observations imply that we may obtain facet-defining inequalities for $\mathcal{P}^2_k$ by \emph{exact lifting}. 

Lifting is a technique for strengthening valid inequalities, and thereby, improving MIP formulations. We illustrate this technique with {an example of lifting EPIs}. Prior to lifting, the facet-defining inequality $w\geq \sum_{i=1}^k \rho_i x_i$ concerns only the variables $x_i$ for $i\in S$. Our goal is to construct the strongest inequality of the form 
\begin{equation}
\label{eq:EPI-lifted}
 w\geq \sum_{i=1}^k \rho_i x_i + \sum_{i=k+1}^n \xi_i x_i. 
\end{equation}
that is valid for $\mathcal{P}^2_k$. To achieve this goal, we go through a \emph{sequence-dependent lifting procedure} to gradually add variables $x_i$, $i\in N\setminus S$, to this inequality with strong coefficients. At each intermediate step, we lift the inequality with $x_j$ for $j\in\{k+1, \dots, n\}$ by solving the $j$-th lifting problem \eqref{eq:EPI_lifting_prob}.
\begingroup
\allowdisplaybreaks
\begin{equation}
\label{eq:EPI_lifting_prob}
\begin{aligned}
\xi_j := \min \hspace{0.2cm} & w - \sum_{i=1}^{k} \rho_i x_i - \sum_{i = k+1}^{j-1} \xi_i x_i&& \\
\textrm{s.t.} \quad & w\geq f\left(a_j + \sum_{i = 1}^{j-1} a_i x_i\right), &&\\
& \sum_{i=1}^{j-1} x_i \leq k-1, && \\
& x \in \{0,1\}^{j-1}. && 
\end{aligned}
\end{equation}
\endgroup
The optimal objective $\xi_j$ is {precisely} the lifted coefficient in the facet-defining inequality $w \geq \sum_{i=1}^k \rho_i x_i + \sum_{i = k+1}^j \xi_i x_i$ for the convex hull of the polyhedron 
\[\mathcal{P}^2_k(\{1,\dots,j\}) = \left\{ (w, x)\in\mathbb{R}\times\{0,1\}^j : w \geq f\left(\sum_{i=1}^j a_ix_i\right), \sum_{i=1}^j x_i \leq k \right\}. \] 
In \cite{yu2023strong}, we give the \emph{closed-form} optimal solutions and objectives of the lifting problems for EPIs and SIs. We propose three classes of strong valid inequalities, which we call the \emph{lifted-EPIs} (LEPIs), \emph{lower-SIs} (LSIs) and \emph{higher-SIs} (HSIs). We refer the readers to \cite{yu2023strong} for the exact forms of these proposed inequalities.

 \begin{example} (LEPIs versus ALIs \cite{yu2023strong})
Suppose $N = [6]$ and $k = 2$. We are given $\mathbf{a} = [ 4, 100,  100,  100,   4,   4]$ and a concave function $f(\mathbf{a}^\top \mathbf{x}) = \sqrt{\mathbf{a}^\top \mathbf{x}}$. The ALI \cite{yu2017polyhedral} with permutation $\bfs{\delta} = (5,2,3,1,4,6)$ is \[w \geq 0.198x_1+ 8.198x_2+ 4.142x_3+ 4.142x_4+ 2x_5+ 0.198x_6. \] The LEPI, lifed from the EPI associated with $S = \{2,5\}$ and $\bfs{\delta}$, is 
\[w\geq 0.828x_1+ 8.198x_2+ 5.944x_3+5.944x_4+ 2x_5+ 0.828x_6.\] 
The lifted-EPI dominates the ALI. In general, lifted-EPIs are at least as strong as the \textit{approximate lifted inequalities} \cite{yu2017polyhedral}. 
\end{example}

\textbf{Cardinality-constrained mean-risk minimization.} We now showcase how the proposed inequalities can be used in a branch-and-cut algorithm to improve computational performance. We first motivate the mean-risk minimization problem with its application in investment portfolio management. Suppose we are given a set, $N$, of $n$ investment options. We have a limited budget which allows us to select at most $k$ options. Let the losses on the investment options be normal random variables $\tilde{\bfs{p}}$ with mean $\bfs{\mu}$ and covariance $Q$. Our goal is to select up to $k$ options for our portfolio, such that with sufficiently high probability, we attain maximal gain, (or equivalently, minimal loss). With this goal, we obtain the following formulation:  
\begin{equation}
\label{eq:VAR}
 \min_{\mathbf{x} \in \{0,1\}^n}\left\{{w} : \mathbb{P}\left( \tilde{\bfs{p}}^\top \mathbf{x} \leq {w} \right) \geq \beta,  \sum_{i=1}^n x_i \leq k\right\}. 
 \end{equation}
Each binary decision variable $x_i$, $i\in N$, represents whether the corresponding option is chosen or not. The auxiliary variable {$w$} captures the total loss of our portfolio, which we would like to minimize. The parameter $\beta\in (0.5, 1)$ represents how risk-averse we are. We may set higher $\beta$ if we prefer lower risk while allowing potentially lower profit, and vice versa. In fact, the optimization problem above is called \emph{Value-at-Risk minimization}, and it can be reformulated as the following \cite{atamturk2008polymatroids,birge2011introduction,atamturk2019lifted,atamturk2020submodularity}:
\begin{equation}
\label{eq:mean_risk}
\min_{\mathbf{x} \in \{0,1\}^n} \left\{-\bfs{\mu}^\top \mathbf{x} + \Omega\sqrt{\mathbf{x}^\top Q\mathbf{x}} : \sum_{i=1}^n x_i \leq k\right\}. 
\end{equation}
Here, $\Omega$ is a constant parameter $\Phi^{-1}(\beta)$, where $\Phi$ is the standard normal cumulative distribution function. We note that $Q$ is a positive semidefinite matrix. Problem \eqref{eq:mean_risk} is \emph{cardinality-constrained mean-risk minimization} with correlated random variables \cite{atamturk2019lifted, atamturk2020submodularity}. We use the notation $\text{diag}(\bfs{y})$ to denote the diagonal matrix with a vector $\bfs{y}$ as its diagonal. The covariance matrix $Q$ is the sum of $Q-\text{diag}(\mathbf{a})$ and $\text{diag}(\mathbf{a})$, for $\mathbf{a}\in\mathbb{R}_+^n$ with $Q-\text{diag}(\mathbf{a}) \succeq 0$. The benefit of this decomposition is that the separable quadratic term $\mathbf{x}^\top \text{diag}(\mathbf{a}) \mathbf{x}$ can be linearized to $\mathbf{a}^\top \mathbf{x}$. As such, problem \eqref{eq:mean_risk} has an equivalent formulation (SOCP): 
\begin{equation}
\label{eq:SOCP}
\begin{aligned}
\min_{(w,y,z,\mathbf{x}) \in \mathbb{R}_+^3\times \{0,1\}^n}  \{-\bfs{\mu}^\top \mathbf{x} + \Omega z : & \;  w \geq \sqrt{\sum_{i\in N} a_i x_i}, \; \sum_{i=1}^n x_i \leq k, \\
& \;  y\geq \sqrt{\mathbf{x}^\top (Q-\text{diag}(\mathbf{a})) \mathbf{x}}, z^2\geq w^2 + y^2\}, 
\end{aligned}
\end{equation}
which can be handled by optimization solvers, such as Gurobi. We notice that the set $\{(w,\mathbf{x}) \in \mathbb{R} \times \{0,1\}^n : w \geq \sqrt{\sum_{i\in N} a_i x_i}, \; \sum_{i=1}^n x_i \leq k\}$, which is defined by a subset of the constraints from Problem \eqref{eq:SOCP}, coincides with $\mathcal{P}_k^m$ given in \eqref{eq:P}. Under this observation, we may use the valid inequalities from \cite{yu2017polyhedral, yu2023strong} in a branch-and-cut algorithm in order to improve the computational performance. To apply LEPIs and LSIs proposed by \cite{yu2023strong}, we decompose $\mathbf{a}$ into $\mathbf{a}^{\text{two}}+\mathbf{a}^{\text{res}}$ so that $\mathbf{a}^{\text{two}}, \mathbf{a}^{\text{res}} \in\mathbb{R}_+^n$, and $\mathbf{a}^{\text{two}}$ contains two distinct weights. Then Problem \eqref{eq:SOCP} can be rewritten as
\begin{align*}
\min_{(v,w,y,z,\mathbf{x}) \in \mathbb{R}_+^4\times \{0,1\}^n} \{-\bfs{\mu}^\top \mathbf{x} + \Omega z : & \: w \geq \sqrt{\sum_{i\in N} {a}^{\text{two}}_i x_i}, \sum_{i=1}^nx_i \leq k, \\
& \: v \geq \sqrt{\sum_{i\in N} {a}^{\text{res}}_i x_i}, y\geq \sqrt{\mathbf{x}^\top (Q-\text{diag}(\mathbf{a})) \mathbf{x}}, \\
&\;  z^2\geq v^2 + w^2 + y^2  \}. 
\end{align*}

Table \ref{table:multi_weight_result} summarizes the computational performance on problem \eqref{eq:mean_risk} of a branch-and-cut algorithm with LEPIs and LSIs (BC-LEPI-LSI), a branch-and-cut algorithm with ALIs (BC-ALI), and directly solving Problem \eqref{eq:SOCP} (SOCP). All test instances are randomly generated. The experiments are executed on one thread of a Linux server with Intel Haswell E5-2680 processor at 2.5GHz and 128GB of RAM. All the solution methods are implemented in Python 3.6 and Gurobi Optimizer 9.5.1. {The first column of Table \ref{table:multi_weight_result} reports the risk tolerance parameter $\beta\in\{0.95, 0.975, 0.99\}$. The second column shows the cardinality upper bound $k\in \{5, 10, 15\}$ in each test case. The fourth column presents the average runtime in seconds. The next column lists the average end gaps. An end gap is computed by (UB$-$LB)/UB$\times$100\%, where UB and LB are the best upper and lower bounds on the objective. The last two columns report the average numbers of branch-and-bound nodes visited and inequalities added.} BC-LEPI-LSI algorithm outperforms BC-ALI and SOCP in all the test cases---BC-LEPI-LSI solves to optimality all but one test case with $\beta = 0.99$ and $k=15$. In this particular case, BC-LEPI-LSI attains a small end gap of 0.6\%. BC-ALI and SOCP have longer average runtime than BC-LEPI-LSI and more failed instances.

\begin{table}[h]
\begin{center}
\small
\caption{\small Computational performance of BC-LEPI-LSI, BC-ALI and SOCP on Problem \eqref{eq:mean_risk}. Superscript $^{i}$ means that $i$ our of five instances are solved within the time limit of one hour. For BC-LEPI-LSI, the average number of total cuts is represented as $\text{m}^{\text{LEPI}} + \text{m}^{\text{LSI}} = \text{m}$, where $\text{m}^{\text{LEPI}}$ is the average number of LEPIs added across five trials, and $\text{m}^{\text{LSI}}$ is that of LSIs.}
\label{table:multi_weight_result}
\begin{tabular}{|c|c|c|c|c|c|c|}
\hline
$\beta$ & $k$ & method & time (s) & end gap & \# nodes & \# cuts \\
\hline
\multirow{9}{*}{0.95}   & \multirow{3}{*}{5}   & BC-LEPI-LSI & $171.9^5$	&	0.0\%	&	6295.2	&	394.6+140.2=534.8\\
                                   &                               & BC-ALI          & $244.1^5$	&	0.0\%	&	10723.2	&	1071.4\\
                                   &                               & SOCP            & $1702.1^5$	&	0.0\%	&	173698.4	&	N/A\\
                                    \cline{2-7}
                                   & \multirow{3}{*}{10} & BC-LEPI-LSI & $577.5^5$	&	0.0\%	&	15612.2	&	1251.0+248.6=1499.6\\
                                   &                               & BC-ALI          & $958.3^5$	&	0.0\%	&	26572.0	&	2656.6\\
                                   &                               & SOCP            & --$^0$	&	5.8\%	&	82144.0	&	N/A\\
                                    \cline{2-7}
                                   & \multirow{3}{*}{15} & BC-LEPI-LSI & $124.3^5$	&	0.0\%	&	3414.6	&	258.0+73.0=331.0\\
                                   &                               & BC-ALI          & $217.2^5$	&	0.0\%	&	5907.4	&	590.0\\
                                   &                               & SOCP            & $2789.6^2$	&	1.7\%	&	79930.6	&	N/A\\
\hline    
\multirow{9}{*}{0.975} & \multirow{3}{*}{5}   & BC-LEPI-LSI & $827.4^5$	&	0.0\%	&	20236.2	&	1176.2+705.2=1881.4\\
                                   &                               & BC-ALI          & $1254.1^5$	&	0.0\%	&	41768.4	&	4175.6\\
                                   &                               & SOCP            & $3349.5^1$	&	31.1\%	&	166471.6	&	N/A\\
                                    \cline{2-7}
                                   & \multirow{3}{*}{10} & BC-LEPI-LSI & $838.9^5$	&	0.0\%	&	21996.8	&	1498.6+648.8=2147.4\\
                                   &                               & BC-ALI          & $1187.4^4$	&	1.2\%	&	40896.2	&	4088.8\\
                                   &                               & SOCP            & --$^0$	&	12.9\%	&	74527.4	&	N/A\\  
                                    \cline{2-7}
                                   & \multirow{3}{*}{15} & BC-LEPI-LSI & $988.6^5$	&	0.0\%	&	19225.8	&	1575.0+324.2=1899.2\\
                                   &                               & BC-ALI          & $1912.7^3$	&	0.8\%	&	43951.8	&	4393.8\\
                                   &                               & SOCP            & --$^0$	&	6.9\%	&	52940.8	&	N/A\\
\hline  
\multirow{9}{*}{0.99}   & \multirow{3}{*}{5}   & BC-LEPI-LSI & $650.1^5$	&	0.0\%	&	20280.0	&	1175.6+715.4=1891.0\\
                                   &                               & BC-ALI          & $947.5^5$	&	0.0\%	&	38703.6	&	3869.4\\
                                   &                               & SOCP            & $3355.9^1$	&	72.3\%	&	259287.8	&	N/A\\
                                    \cline{2-7}
                                   & \multirow{3}{*}{10} & BC-LEPI-LSI & $1855.7^5$	&	0.0\%	&	38376.0	&	2946.4+842.4=3788.8\\
                                   &                               & BC-ALI          & $3068.4^1$	&	4.7\%	&	88635.6	&	8862.6\\
                                   &                               & SOCP            & --$^0$	&	30.8\%	&	75435.8	&	N/A\\
                                    \cline{2-7}
                                   & \multirow{3}{*}{15} & BC-LEPI-LSI & $1354.1^4$	&	0.6\%	&	24718.6	&	2269.6+184.4=2454.0\\
                                   &                               & BC-ALI          & $2434.1^2$	&	1.8\%	&	70845.8	&	7083.4\\
                                   &                               & SOCP            & --$^0$	&	10.5\%	&	65965.4	&	N/A\\
\hline                        
\end{tabular}
\end{center}
\end{table}

The previous sections showcase the power of submodular set function optimization in the modeling and computation of decision-making problems. Next, we elaborate on two notions of generalized submodularity, namely the multi-set extension and the mixed-integer extension, and how the techniques described previously can be applied to these cases.

\section{Multi-Set Extension of Submodularity} 
\label{sect:multiset}
In this section, we consider the situations in which we are given heterogenous items, and our task is to determine the best subset of the $\lq$bins' to fill, as well as the types of items to place in the selected $\lq$bins'  (see Figure \ref{fig:sub_extensions}). In particular, we discuss the multi-set extension of submodularity, called \emph{$k$-submodularity}. Here, $k\geq 1$ is a positive integer that represents the number of item types. 

\begin{figure}[h] 
   \centering
    \includegraphics[width=6cm]{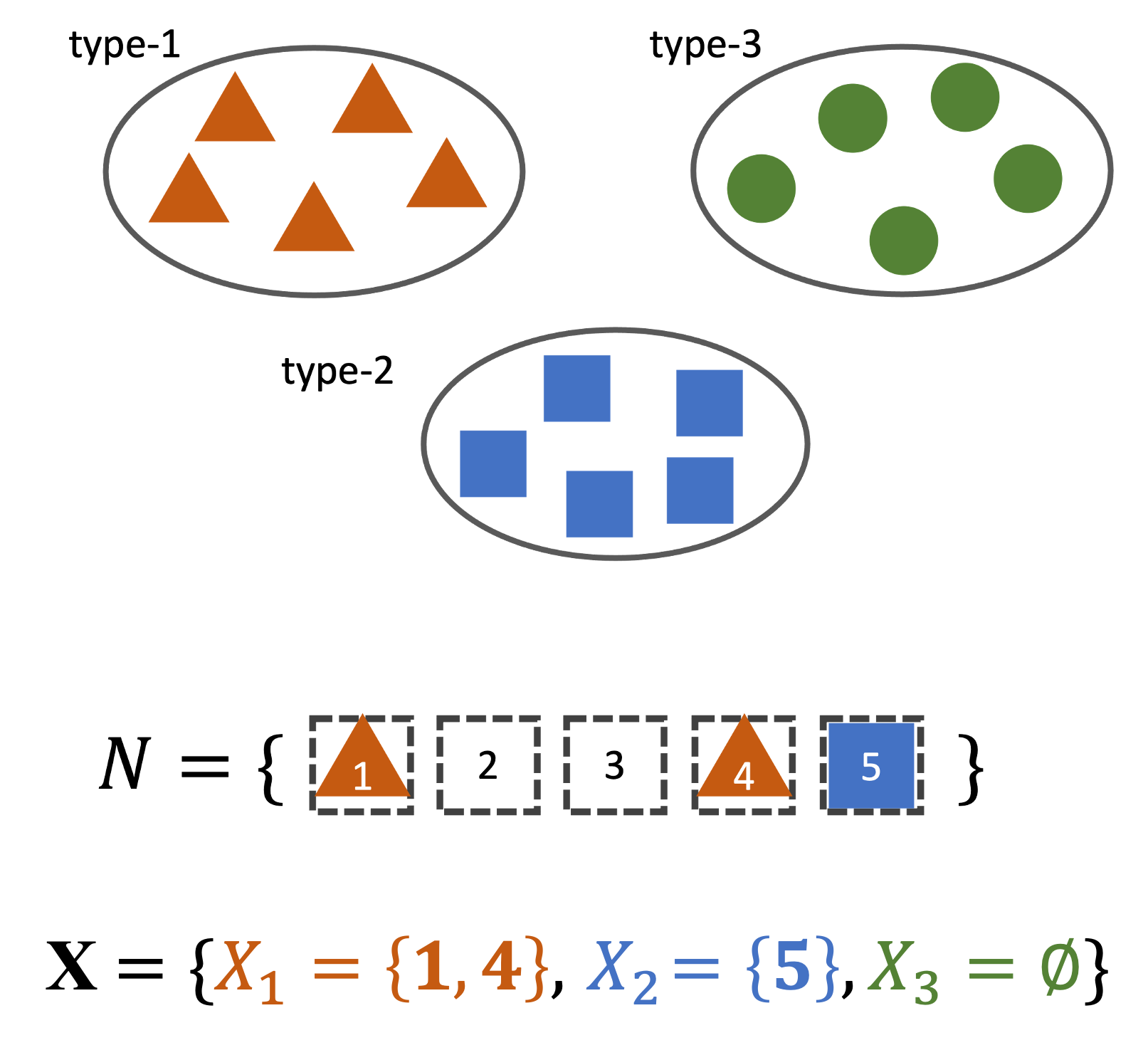}
    \caption{An example of a $k$-set. }
    \label{fig:k-set}
\end{figure}

Before formally defining $k$-submodularity, we first set forth some notations. Recall that $N = \{1,2,\dots, n\}$ is a finite and non-empty ground set. We call a $k$-tuple of pairwise disjoint subsets of $N$ a \emph{$k$-set}. Figure \ref{fig:k-set} provides an example of a $3$-set---we are given three types of items and $\mathbf{X}$ is a $3$-set. Each component $X_i$, $i\in\{1,2,3\}$, of $\mathbf{X}$ contains the $\lq$bins' at which items of type $i$ are placed. Every $\lq$bin' can hold one item, so these components are pairwise disjoint. We use  
\[(k+1)^N = \{(S_1,S_2, \dots, S_k): S_q\subseteq N, S_q\cap S_{q'} = \emptyset, \text{ for all } q, q' \in \{1,2,\dots,k\} \text{ with } q \neq q'\}\] 
to denote the collection of all $k$-sets. For ease of notation, we let $\mathbf{S} = (S_1,S_2, \dots, S_k)\in (k+1)^N$. When $k=1$, $(k+1)^N = 2^N$ is exactly the power set of $N$.

\begin{definition}
\label{def:k_sub}
For any positive integer $k$, a function $f:(k+1)^N \rightarrow \mathbb{R}$ is \emph{$k$-submodular} if for any $\mathbf{X}=(X_1, X_2, \dots, X_k), \mathbf{Y}=(Y_1, Y_2, \dots, Y_k)\in (k+1)^N$, 
\begin{align*}
f(\mathbf{X}) + f(\mathbf{Y}) \geq &\hspace{0.2cm} f(\mathbf{X} \sqcap \mathbf{Y}) + f(\mathbf{X} \sqcup \mathbf{Y}),
\end{align*}
where \[\mathbf{X} \sqcap \mathbf{Y} = (X_1\cap Y_1, X_2 \cap Y_2, \dots, X_k \cap Y_k),\] and \[\mathbf{X} \sqcup \mathbf{Y} = \left((X_1\cup Y_1 )\backslash \bigcup_{q=2}^k(X_q\cup Y_q), \dots, (X_k \cup Y_k)\backslash \bigcup_{q=1}^{k-1} (X_q\cup Y_q)\right).\] 
\end{definition} 

Definition \ref{def:k_sub} is similar to Definition \ref{def:submodular} for submodular set functions. They differ by the fact that the intersection and union operators in Definition \ref{def:k_sub} apply to $k$-sets component-wise. In particular, the $\sqcup$ operator excludes additional items to ensure pairwise disjointness of the resulting $k$-set. When $k=1$, the functions satisfying Definition \ref{def:k_sub} are exactly the submodular set functions. When $k=2$, we call these functions \emph{bisubmodular functions}.

\subsection{Applications of $k$-submodularity}
Next, we describe a few applications of $k$-submodularity to demonstrate its modeling power. 

\subsubsection{Multi-type sensor placement.} 
\label{subsect:sensor}
Internet of things (IoT) has enabled sensor networks to monitor systems such as smart homes \cite{ghayvat2015wellness}, smart grids \cite{abujubbeh2019software} and smart cities \cite{zanella2014internet}. In these applications, we often have more than one type of sensors that measure different aspects of the environment (e.g. temperature, humidity, and luminance). Suppose we are given $k\geq 2$ types of sensors and a set $N = \{1,2,\dots,n\}$ of candidate locations such that each location holds up to one sensor. Let the locations for type-$i$ sensors be $S_i$. Then every $k$-set, $\mathbf{S} = (S_1, \dots, S_k)\in (k+1)^N$, uniquely represents a $k$-type sensor placement plan. The effectiveness of any placement plan can be quantified using the following entropy function
$H(X_{\mathbf{S}}) = -\sum_{x\in \mathcal{X}_{\mathbf{S}}} \mathbb{P}(X_{\mathbf{S}} = x)\log \mathbb{P}(X_{\mathbf{S}} = x)$, 
where $X_{\mathbf{S}}$ is a random variable representing the possible observations made by sensors of corresponding types installed at $\mathbf{S}$, and the set $\mathcal{X}_{\mathbf{S}}$ contains all possible observations \cite{ohsaka2015monotone}. Intuitively, the entropy of $X_{\mathbf{S}}$ is high when the $k$ types of measurements at locations $S_1,\dots, S_k$, respectively, are difficult to predict. The function $f(\mathbf{S}) = H(X_{\mathbf{S}})$ for all $\mathbf{S}\in (k+1)^N$, is $k$-submodular \cite{ohsaka2015monotone}.

\subsubsection{Multi-topic influence propagation.}
Social media platforms have served as a crucial source of information. The mechanism of information spread over a social network can be modeled using a submodular set function, which maps a set of initial information spreaders to the expected total number of network users who are influenced by this information \cite{kempe2015maximizing}. This model is particularly helpful for applications such as viral marketing. This model is generalized to allow $k\geq 2$ types of influence \cite{ohsaka2015monotone}. In a social network $G = (N, A)$, $N=\{1,2,\dots, n\}$ represents the network users, and the arc set $A$ models the interactions among the users. For any $i,j\in N$ and $q\in\{1,2,\dots,k\}$, $p_{(i,j)}^q$ is the probability of $i$ influencing $j$ on topic $q$. In viral marketing, $p_{(i,j)}^q$ can be interpreted as the likelihood of $i$ convincing $j$ to purchase product $q$. Let $S_q\subseteq N$ be the influencers selected to promote the type-$q$ products. Moreover, let $A_q(S_q)$ be the individuals influenced by the initial spreaders $S_q$ about product $q$ under the stochastic model proposed by \cite{kempe2015maximizing}. The function $f(\mathbf{S}) = \mathbb{E}[ |\bigcup_{q=1}^k A_q(S_q)| ]$ evaluates the expected total number of influenced individuals given initial influencers $\mathbf{S}$, and it has been shown to be monotone and $k$-submodular \cite{ohsaka2015monotone}.

\subsubsection{Multi-class feature selection.} 
Feature selection is an important technique in numerous fields of research, such as machine learning \cite{shalev2014understanding}, bioinformatics \cite{saeys2007review}, and data mining \cite{rokach2008data}. This process improves the interpretability of the models and avoids the curse of dimensionality. Multi-class feature selection generalizes the usual feature selection problem. In $k$-class feature selection, we are given $k$ uncorrelated prediction variables and a pool of features. Our goal is to select the most informative features and to determine the corresponding prediction variables for the selected features. {When $k=2$, the problem is called the coupled feature selection problem \cite{singh2012bisubmodular}}. More specifically, $C_1, C_2$ are two prediction variables and $N$ is a set of features. For any $\mathbf{S} = (S_1, S_2)\in 3^N$, $S_i$ contains the features selected for predicting $C_i$, $i \in \{1,2\}$. It is assumed that $S_1, S_2$ are mutually conditionally independent given $C = \{C_1,C_2\}$. The informativeness of a coupled feature selection is measured by its mutual information 
$I(\mathbf{S}; C) = H(S_1\cup S_2) - \sum_{i\in S_1}H(i\mid C_1) - \sum_{j\in S_2}H(j\mid C_2)$, 
where $H(X) = -\sum_{x\in \mathcal{X}} \mathbb{P}(X = x)\log \mathbb{P}(X = x)$, and $H(X\mid Y) = -\sum_{x\in \mathcal{X},y\in\mathcal{Y}} \mathbb{P}(X=x, Y=y) \log \frac{\mathbb{P}(X=x, Y=y)}{\mathbb{P}(Y = y)}$. Higher mutual information suggests that the biset of features is more informative for the prediction tasks. The function $f(\mathbf{S}) = I (\mathbf{S};C)$ is monotone and bisubmodular \cite{singh2012bisubmodular}. 

\subsubsection{Drug-drug interaction detection.}
When certain drugs are taken concomitantly, they react in our bodies. Such reactions are referred to as the drug-drug interactions (DDIs). DDIs comprise at least 30\% of the adverse drug events, which cause 770,000 injuries and deaths every year \cite{tatonetti2012novel, pirmohamed1998drug}. {The authors in} \cite{hu2019bi} introduce a function that models the correlation between combinations of drugs and adverse effects, and they further show that this function is bisubmodular and determine DDIs based on online medical forum data by exploiting bisubmodularity.

\subsection{Exact solution methods for $k$-submodular optimization}
The previous section illustrates how $k$-submodular functions model utilities in real-world applications. After obtaining such a model, our next step is to optimize the $k$-submodular utility function in order to determine an optimal decision. In this section, we describe exact solution approaches to $k$-submodular function maximization and minimization. Throughout this section, we assume without loss of generality that $f(\bfs{\emptyset}) = 0$, where $\bfs{\emptyset}=(\emptyset, \dots, \emptyset)$.

\subsubsection{$k$-submodular maximization.}
A $k$-submodular maximization problem assumes the following form: 
\begin{equation}
\label{eq:original_max}
\max_{\mathbf{X} \in \mathcal{X}} f(\mathbf{X}),
\end{equation}
where $f$ is $k$-submodular, and $\mathcal{X}\subseteq (k+1)^N$ contains all the feasible $k$-sets. When Problem \eqref{eq:original_max} is unconstrained, $\mathcal{X} = (k+1)^N$. We may also incorporate linear constraints, such as cardinality and knapsack constraints, to make $\mathcal{X}\subsetneq (k+1)^N$. 

Submodular set function maximization is a special case of $k$-submodular maximization. Given that submodular maximization is NP-hard, $k$-submodular maximization is NP-hard as well. Most of the methods proposed for $k$-submodular maximization are approximation algorithms. {A class of bisubmodular functions can be maximized by a constant-factor approximation algorithm \cite{singh2012bisubmodular}. Studies including} \cite{iwata2013bisubmodular, ward2014maximizing, ward2016maximizing} establish a $1/2$-approximation guarantee of a randomized greedy algorithm for unconstrained bisubmodular maximization. When $k\geq 3$, \cite{ward2016maximizing} achieve an approximation guarantee of $\max(1/3, 1/(1+\max(1, \sqrt{(k-1)/4})))$, and this result has been improved to $1/2$ by  \cite{iwata2016improved}. Other studies consider maximizing nonnegative monotone $k$-submodular functions under various constraints (e.g., cardinality constraints \cite{ohsaka2015monotone}, matroid constraints \cite{sakaue2017maximizing}). 

In this tutorial, we focus on the exact solution methods derived from mixed-integer programming approaches. By slightly abusing notation, we introduce the following binary representation of the $k$-sets. That is, every $\mathbf{X}\in(k+1)^N$ is uniquely mapped to $\mathbf{x} = [{\mathbf{x}^1}, \dots, {\mathbf{x}^k}]^\top\in \{0,1\}^{kn}$, with $\sum_{q = 1}^k x^q_i \leq 1 \text{ for all } i\in N$ and 
\[x^q_i =\begin{cases}
1, & \text{if } i\in X_q, \\
0, & \text{otherwise.} 
\end{cases}\] 
The linear inequalities $\sum_{q = 1}^k x^q_i \leq 1$ for all $i\in N$ ensure pairwise disjointness of the components in the corresponding $k$-set. We let $\mathcal{X}$ denote the feasible $k$-sets, as well as the feasible binary incidence vectors. Then we can rewrite Problem \eqref{eq:original_max} as 
\begin{equation}
\label{eq:hypograph}
\max\{w : (\mathbf{x}, w)\in \mathcal{T}_f, \mathbf{x}\in \mathcal{X}\}, 
\end{equation} 
where \[\mathcal{T}_f = \left\{(\mathbf{x}, w)\in \{0,1\}^{kn} \times \mathbb{R} : w \leq f(\mathbf{x}), \sum_{q = 1}^k x^q_i \leq 1 \text{ for all } i\in N \right\}. \]

Let $\rho_{q,i}(\mathbf{X}) = f(X_1,\dots, X_q \cup \{i\}, \dots, X_k) - f(\mathbf{X})$ be the marginal contribution of adding item $i\in N$ to the $q$-th component of $\mathbf{X}\in (k+1)^N$. 

\begin{definition} 
A $k$-submodular function $f: (k+1)^N\rightarrow \mathbb{R}$ is \emph{monotone non-decreasing}, or simply \emph{monotone}, if $f(\mathbf{Y})\geq f(\mathbf{X})$ for any $\mathbf{X},\mathbf{Y}\in (k+1)^N$ such that $X_q \subseteq Y_q$ for all $q \in \{1,\dots, k\}$. Equivalently, $f$ is monotone if $\rho_{q,i}(\mathbf{X}) \geq 0$, for any $\mathbf{X} \in (k+1)^N$, $i\in N\backslash \bigcup_{p=1}^k X_p$, and $q \in \{1,\dots, k\}$.
\end{definition}

Given a non-monotone $k$-submodular function $f$, we can construct its monotone counterpart as the following. For every $i\in N$ and $q\in\{1,\dots,k\}$, let
\begin{equation*}
\label{eq:xi}
\xi_i^q = \min\left\{\rho_{q, i}(\mathbf{S}) : \mathbf{S}\in (k+1)^{N\backslash\{i\}}, \bigcup_{p=1}^k S_p = N\backslash\{i\} \right\}.
\end{equation*}

\begin{lemma} \cite{yu2021exact}
Given any non-monotone $k$-submodular function $f$, 
\begin{align*}
f^*(\mathbf{X}) := f(\mathbf{X}) - \sum_{q=1}^k \sum_{i\in X_q} \xi_i^q
\end{align*}
is $k$-submodular and monotone. 
\end{lemma}
We refer the readers to Lemma 3.2 of \cite{yu2021exact} for more details on this result. 

\begin{definition}
Given any $k$-submodular function $f$ and any $\mathbf{S}\in (k+1)^N$, the \emph{$k$-submodular inequality} associated with $\mathbf{S}$ is 
\begin{equation}
\label{cut_general}
\begin{aligned}
w \leq f(\mathbf{S})  & + \sum_{q=1}^k \sum_{i\notin \bigcup_{r=1}^k S_r} \rho_{q, i}(\mathbf{S})x_i^q  +  \sum_{q=1}^k \sum_{p\in\{1,\dots, k\}\backslash \{q\}} \sum_{i\in S_p} \rho_{q,i}(\bfs{\emptyset})x_i^q - \sum_{q=1}^k \sum_{i\in S_q} \xi_i^q (1-x_i^q). 
\end{aligned}
\end{equation}
\end{definition} 

\begin{proposition} \cite{yu2021exact}
For any $k$-submodular function $f$ and any $\mathbf{S}\in (k+1)^N$, the $k$-submodular inequality \eqref{cut_general} associated with $\mathbf{S}$ is valid for $\mathcal{T}_f$. 
\end{proposition}
We refer the readers to Proposition 4.2 in \cite{yu2021exact} for the proof of this proposition. We note that the submodular inequalities \eqref{eq:sub_ineq} proposed by \cite{wolsey1999integer} is a special case of the $k$-submodular inequalities when $k=1$. The $k$-submodular inequalities provide a piecewise linear representation of the objective function $f$, with which we derive another equivalent formulation of $k$-submodular maximization as shown in the next theorem. 

\begin{theorem} \cite{yu2021exact}
Problem \eqref{eq:hypograph} is equivalent to 
\begin{equation}
\label{eq:k_sub_form}
\begin{aligned}
\max_{(\mathbf{x}, w)\in \{0,1\}^{kn}\times \mathbb{R}} \{w : & \; w \leq  f(\mathbf{S}) + \sum_{q=1}^k \sum_{i\notin \bigcup_{r=1}^k S_r} \rho_{q, i}(\mathbf{S})x_i^q  +  \sum_{q=1}^k \sum_{p\in\{1,\dots, k\}\backslash \{q\}} \sum_{i\in S_p} \rho_{q,i}(\bfs{\emptyset})x_i^q \\
 & - \sum_{q=1}^k \sum_{i\in S_q} \xi_i^q (1-x_i^q),  \forall \;\mathbf{S}\in (k+1)^N, \\
 & \sum_{q=1}^k x^q_i \leq 1, \forall \; i\in N, \mathbf{x} \in \mathcal{X}.\}
 \end{aligned}
 \end{equation}
\end{theorem}

We use the DCG framework (see Algorithm \ref{alg:general_dcg}) to overcome the challenge of having exponentially many $k$-submodular inequalities in the mixed 0-1 linear program reformulation \eqref{eq:k_sub_form}. A computational study on multi-type sensor placement problem is included in Section \ref{sect:sensor}.

\subsubsection{$k$-submodular minimization}
In this section, we summarize the results on \emph{minimizing $k$-submodular functions}. In particular, we focus on the case where $k=2$. {The study} \cite{qi1988directed} proves an analogue of Lov\'asz extension (see Section \ref{sect:preliminaries}) for bisubmodular functions. This result suggests that unconstrained bisubmodular minimization can be solved in polynomial time using the ellipsoid method. Weakly and strongly polynomial algorithms have also been proposed for unconstrained bisubmodular minimization \cite{fujishige2005bisubmodular,mccormick2010strongly}. The Min-Max Theorem for submodular and bisubmodular minimization is generalized to the $k$-submodular case with $k\geq 3$ \cite{huber2012towards}, but the complexity of general $k$-submodular minimization remains an open problem. In \cite{yu2020polyhedral}, we completely characterize the epigraph convex hull of any bisubmodular function, and this polyhedral characterization leads to an effective cutting plane algorithm to solve \emph{constrained} bisubmodular minimization to global optimality. 

We now consider  the constrained bisubmodular minimization problem:  
\begin{equation}
\label{eq:original_bisub}
\begin{aligned}
\min & \quad f(S_1,S_2) \\
\text{s.t.} & \quad (S_1,S_2) \in \mathcal{S} \subseteq 3^N, 
\end{aligned}
\end{equation}
and summarize the mixed-integer programming approach {we propose in} \cite{yu2020polyhedral}. Without loss of generality, we again assume $f(\emptyset,\emptyset) = 0$. In addition, $\mathcal{S}\subseteq 3^N$ is the set of feasible solutions after possibly incorporating any linearly representable constraints. This is a class of nonlinear optimization problems over bisets, which could not be directly entered into optimization solvers.  

In order to reformulate Problem \eqref{eq:original_bisub} as a mixed-integer linear program, we first convert the bisets into ternary characteristic vectors. To be more precise, for any $(S_1,S_2) \in 3^N$, there is a unique vector $\mathbf{x}\in \{0,\pm 1\}^n$ such that for every $i\in \{1,2,\dots, n\}$, 
\[ x_i := \begin{cases}
1, & i\in S_1, \\
-1, & i\in S_2, \\
0, & \text{otherwise.}
\end{cases}\]
Conversely, every $\mathbf{x}\in \{0,\pm1\}^n$ corresponds to a unique biset. To this point, we have successfully turned the set $\mathcal{S}$ of feasible bisets in Problem \eqref{eq:original_bisub} to the integer vectors $\mathbf{x}\in \mathcal{X}$, where $\mathcal{X}$ contains every $\mathbf{x}\in \{0,\pm 1\}^n$ that also satisfies all the additional linear constraints. 

The next step is to rewrite the nonlinear objective function $f$ into a form that the optimization solvers could handle. In \cite{yu2020polyhedral}, we propose a {piecewise linear} representation of $f$. This representation relies on the \emph{extremal poly-bimatroid inequalities}, which are defined as the following. 

\begin{definition}
Let any permutation $\bfs{\delta} = (\delta_1, \delta_2, \dots, \delta_n)$ of $N$ and any $\bfs{\sigma}\in\{\pm 1\}^n$ be given. 
An \emph{extremal poly-bimatroid inequality} associated with $\bfs{\delta}$ and $\bfs{\sigma}$ is $w \geq \bfs{\pi}^\top \mathbf{x}$, where $\bfs{\pi} \leftarrow \emph{\texttt{Signed\_Greedy}}(\bfs{\delta}, \bfs{\sigma})$ (see Algorithm \ref{alg:signed_greedy}).  
\end{definition}

\begin{algorithm}[H]
 \caption{\texttt{Signed\_Greedy} \cite{ando1996structures, mccormick2010strongly}}
 \label{alg:signed_greedy} 
\begin{algorithmic}[1]
\STATE \textbf{Input} a permutation $\bfs{\delta} = \{\delta_1, \delta_2, \dots, \delta_n\}$ of $N$, a vector $\bfs{\sigma}\in\{\pm 1\}^n$\;
\STATE $\bfs{\pi}  \leftarrow \mathbf{0}$, $S_1 \leftarrow \emptyset$, $S_2 \leftarrow \emptyset$\;
\FOR {$i=1,2,\dots,n$}
    \IF {$\sigma_{\delta_i} = 1$}
        \STATE $\pi_{\delta_i} \leftarrow f(S_1\cup\{\delta_i\}, S_2) - f(S_1, S_2)$\;
        \STATE $S_1 \leftarrow S_1\cup \{\delta_i\}$\;
    \ELSE
        \STATE  $\pi_{\delta_i} \leftarrow -f(S_1, S_2\cup\{\delta_i\}) + f(S_1, S_2)$\;
        \STATE $S_2 \leftarrow S_2\cup \{\delta_i\}$\;
    \ENDIF
\ENDFOR
\STATE  \textbf{Output} $\bfs{\pi} \in\mathbb{R}^n$. 
\end{algorithmic}
\end{algorithm}

\vspace{0.4cm} Theorem \ref{thm:bisub_conv} below formally establishes that the extremal poly-bimatroid inequalities provide a {piecewise linear} representation of $f$. We refer the readers to \cite{yu2020polyhedral} for the proof of Theorem \ref{thm:bisub_conv}. 
\begin{theorem}
\label{thm:bisub_conv}
\cite{yu2020polyhedral} The convex hull of the epigraph of $f$, $\conv{(\mathbf{x} , w)\in \{0, \pm 1\}^n\times \mathbb{R} : w \geq f(\mathbf{x})}$, is described completely by all the extremal poly-bimatroid inequalities $w \geq \bfs{\pi}^\top \mathbf{x}$, and the trivial inequalities $-1\leq x_i \leq 1$ for all $i\in \{1,2,\dots, n\}$. 
\end{theorem}

Following from Theorem \ref{thm:bisub_conv}, the constrained bisubmodular minimization problem \eqref{eq:original_bisub} has an equivalent reformulation: 
\begin{equation}
\label{eq:reform_bisub}
\begin{aligned}
\min & \quad w  \\
\text{s.t.} & \quad (\mathbf{x},w)\in \mathcal{C}, \\
& \quad \mathbf{x} \in \mathcal{X}.
\end{aligned}
\end{equation} 
The polyhedral set $\mathcal{C}$ is defined by all the extremal poly-bimatroid inequalities. The set $\mathcal{X}$ contains all the ternary characteristic vectors that correspond to feasible bisets. Due to the observation that Problem \eqref{eq:reform_bisub} has exponentially many extremal poly-bimatroid inequalities, it is desirable to apply the DCG framework (see Algorithm \ref{alg:general_dcg}). The problem of separating an extremal poly-bimatroid inequaliy at any fractional solution $\overline{\mathbf{x}}\in \mathbb{R}^n$ can be efficiently solved by an $\mathcal{O}(n\log n)$ greedy algorithm (Algorithm \ref{alg:bisub_gen_greedy}). 

\vspace{0.2cm}\begin{algorithm}[H]
\caption{\texttt{Generalized\_Greedy} \cite{bouchet1987greedy}}
\label{alg:bisub_gen_greedy}
\begin{algorithmic}[1]
\STATE \textbf{Input} $\overline{\mathbf{x}}\in \mathbb{R}^n$\;
\STATE Sort entries in $\overline{\mathbf{x}}$ to obtain a permutation $\bfs{\delta}$ such that $|\overline{x}_{\delta_1}| \geq |\overline{x}_{\delta_2}| \geq \dots \geq |\overline{x}_{\delta_n}|$\;
\STATE $\bfs{\pi} \leftarrow \mathbf{0}$, $S_1, S_2 \leftarrow \emptyset$\;
\FOR {$i = 1, 2, \dots, n$}
    \IF {$\overline{x}_{\delta_i}\geq 0$}
        \STATE $\pi_{\delta_i} \leftarrow f(S_1\cup\{\delta_i\}, S_2) - f(S_1, S_2)$\;
        \STATE $S_1 \leftarrow S_1\cup \{\delta_i\}$\;
    \ELSE
        \STATE $\pi_{\delta_i} \leftarrow -  f(S_1, S_2\cup\{\delta_i\}) + f(S_1, S_2)$\;
        \STATE $S_2 \leftarrow S_2\cup \{\delta_i\}$\;
    \ENDIF
\ENDFOR
\STATE  \textbf{Output} An extremal poly-bimatroid inequality $w \geq \bfs{\pi}^\top \mathbf{x}$.
\end{algorithmic}
\end{algorithm}

\vspace{0.4cm}Now that we have described exact solution methods to tackle $k$-submodular maximization and minimization, we demonstrate the proposed DCG algorithms in the multi-type sensor placement problem, as well as its robust variant.

\subsection{An In-Depth Example of Multi-Type Sensor Placement} 
\label{sect:sensor}
We delve deeper into the multi-type sensor placement problem, described in Section \ref{subsect:sensor}, to demonstrate how $k$-submodular optimization is applied in practice. We refer the readers to Example 5.1 in \cite{yu2020polyhedral} for a numerical example of the entropy function {that models the utility in this application}. Recall that we are given $k$ types of sensors and $n$ candidate locations. We would like to determine an optimal $k$-type sensor placement plan $\bfs{S}\in (k+1)^N$ that maximizes entropy. In addition, we have a limited budget which restricts the number of sensors available to us---we may install up to $B_q$ sensors of type $q$, for $q\in\{1,2,\dots,k\}$. These restrictions are formulated as the cardinality constraints, $\sum_{i=1}^n x^q_i \leq B_q$ for every $q\in \{1,2,\dots, k\}$. 

Because of the high nonlinearity of the objective function, there is no compact mixed-integer linear formulation for this problem. Thus we compare the DCG algorithm with $k$-submodular inequalities against the only exact solution approach, exhaustive search (ES). We randomly generate test instances using the Intel Berkeley research lab dataset \cite{bodik2004intel}, which contains the sensor readings of luminance, temperature, and humidity at 54 locations in the lab. For the set of experiments with $k=2$, we aim to find the best placement plan for luminance and temperature sensors. When $k=3$, our goal is to determine the optimal placement plan for all three measurements. Our experiments are executed on two threads of a Linux server with Intel Haswell E5-2680 processor at 2.5GHz and 128GB of RAM. Our algorithms are implemented in Python 3.6 and Gurobi Optimizer 7.5.1 with time limit of one hour. We include our computational results in Table \ref{res:entropy}. We alter $n$, the number of candidate locations, between 20 and 50. We use 100 tuples of luminance, temperature and humidity readings to evaluate entropy. The cardinality upper bounds $B_q$ are set to $n/10$ for all $q\in\{1,\dots, k\}$. We report the computational statistics of DCG, including runtime, the number of $k$-submodular inequalities added, the number of visited branch-and-bound nodes, and the end gaps (i.e., (UB$-$LB)/UB). For ES, we report the average runtime. 

\begin{table}[htb]
\begin{center}
\caption{\small Computational performance of DCG and ES in $k$-type sensor placement. The statistics are averaged across three trials. The superscript $^\ell$ means that $\ell$ instances reach the time limit of one hour.} 
\label{res:entropy}
\begin{tabular}{c|c||c|c|c|c||c}
\hline
 $k$ &  $n$  &      time (s)   &     \# cuts  &     \# nodes & end gap &  ES time (s)  \\ 
  \hline
 \multirow{4}{*}{2} & 20 &	1.09	&	38.67	&	42.33	& -- &	17.68 \\
& 30	&	13.04	&	250.33	&	253.67	& -- &	--$^3$ \\
& 40	&	161.80	&	1783.33	&	1791.00	& -- & 	--$^3$ \\
& 50	&	478.09	&	3559.33	&	3566.33	& -- & 	--$^3$ \\
   \hline
  \multirow{4}{*}{3}		&	20	&	1.58	&	32.67	&	35.33	&	--	&	2953.18 \\
& 	30	&	20.27	&	222.33	&	226.67	&	--	&	--$^3$ \\
	& 40	&	164.97	&	1043.00	&	1049.33	&	--	&	--$^3$ \\
& 50 	&	1779.49	&	6907.67	&	6920.00	&	--	&	--$^3$ \\
\hline
\end{tabular}
\end{center}
\end{table}

DCG solves all test instances to optimality under half an hour. On the other hand, ES fails to solve any instances beyond $n=20$ for both $k=2$ and $k=3$. Consider the test case of $k=2$ and $n=50$ as an example, ES needs to enumerate $50!/(5!5!40!) \approx 2.59\times 10^{12}$ feasible solutions to find an optimal solution. Each entropy evaluation takes $1.6\times 10^{-4}$ seconds on average, so it takes approximately 13.13 years for ES to reach optimality. In contrast, our algorithm finds an optimal solution under eight minutes.

We further consider a robust variant of the coupled sensor placement problem. This time, our goal is to determine a coupled sensor placement plan with the highest worst-case entropy. In particular, we would like our placement plan to be robust against two kinds of potential faults. First, given any placement plan, up to $W\in\mathbb{Z}_+$ sensors of the wrong type will be installed. Second, some sensors will malfunction due to software or hardware faults, but at least $B'_1$ type-1 sensors and $B'_2$ type-2 sensors will function properly. Problem \eqref{eq:robust-entropy} is the robust coupled sensor placement problem. 
\begingroup
\allowdisplaybreaks
\begin{subequations}
\label{eq:robust-entropy}
\begin{alignat}{2}
\max_{(S_1,S_2)\in 3^N} \min_{(T_1, T_2) \in 3^{S_1\cup S_2}} \hspace{0.2cm} & f(T_1, T_2)&& \\
\textrm{s.t.} \quad & |S_1| = B_1,  |S_2| = B_2, &&\\
& |T_1| \geq B'_1,  |T_2| \geq B'_2, && \\
& |T_1\cap S_2| + |T_2\cap S_1| \leq W.  && 
\end{alignat}
\end{subequations}
\endgroup
Given a placement plan $(S_1,S_2)\in 3^N$ selected in the outer maximization problem,  in the inner minimization problem, the decision $(T_1, T_2)$ corresponds to the functioning sensors after the two kinds of faults in the worst case. This inner minimization problem is constrained bisubmodular minimization. Due to the nonlinear objective and cardinality constraints, we apply the DCG algorithm to tackle the mixed-integer reformulation \eqref{eq:bisub_robust} of the inner problem. In this formulation, $\mathcal{C}$ in constraint \eqref{subeq:linearapprox} is defined by all the extremal poly-bimatroid inequalities. Classical robust optimization techniques such as the L-shaped method \cite{laporte1993integer} could be used to tackle the outer maximization problem. 

\begingroup
\allowdisplaybreaks
\begin{subequations}
\label{eq:bisub_robust}
\begin{alignat}{2}
\min \hspace{0.2cm} & w  \\
\textrm{s.t.}\hspace{0.2cm}  & (\mathbf{x} , w) \in \mathcal{C}, && \label{subeq:linearapprox} \\
& x_i = y_i^1 -y_i^2, && \quad \text{ for all } i\in N, \label{subeq:biset}\\
& y_i^1 + y_i^2 \leq 1, &&\quad \text{ for all } i \in N, \label{subeq:disjoint}\\
& \sum_{i\in N} y_i^1 \geq B'_1, && \label{subeq:card1}\\
& \sum_{i\in N} y_i^2 \geq B'_2, && \label{subeq:card2}\\
& y_i^1, y_i^2 = 0, &&\quad \text{ for all } i\in N\backslash (S_1\cup S_2), \label{subeq:outer}\\
& \sum_{i\in S_1}y_i^2 + \sum_{i\in S_2}y_i^1 \leq W, &&\quad \text{ for all } i \in N, \label{subeq:switch}\\
& y_i^1, y_i^2 \in \{0,1\}, &&\quad  \text{ for all } i\in N. 
\end{alignat}
\end{subequations}
\endgroup

We again generate random instances using the temperature and humidity readings from Intel Berkeley research lab dataset. We randomly pick $n$ candidate locations, and we vary the number $t$ of temperature and humidity readings that we use for entropy evaluation. We set $B_1 = \floor{2n/5}$, $B_2 = \floor{n/2}$, $B'_1= \floor{4B_1/5}$, $B'_2=\floor{3B_2/5}$, and $W=\floor{3(B_1+B_2)/5}$. In Table \ref{res:entropy_min}, we report the average runtime of DCG in seconds, the average numbers of extremal poly-bimatroid inequalities added, and the average numbers of branch-and-bound nodes visited. We observe that DCG solves all test instances within five minutes on average. In smaller test cases with $n \leq 10$, DCG attains optimality within five seconds. Overall, DCG is an efficient method to exactly optimize the inner constrained bisubmodular minimization problem. 

\begin{table}[htb]
\small
\begin{center}
\caption{\small Computational results of DCG for the bisubmodular minimization subproblem in robust coupled sensor placement. The statistics are averaged across ten instances.}
\begin{tabular}{c|c|c|c|c}
\hline
$n$ & $t$ & time (s) &   \# cuts & \# nodes \\
\hline
 \multirow{4}{*}{5} & 10 & 0.003  & 2.1  & 2.0 \\
 & 20 & 0.012  & 10.2  & 25.9 \\
 & 50 & 0.011  & 4.9  & 6.4 \\
 & 100 & 0.035  & 11.4  & 28.9 \\
 & 500 & 0.096  &  8.0 & 20.7 \\
\hline
 \multirow{4}{*}{10} & 10 & 0.104  & 62.4  & 493.1\\
 & 20 &  0.238 & 97.1  & 1113.1 \\
 & 50 & 0.339  & 77.9  & 750.0 \\
 & 100 & 0.896  &  114.9 & 1024.1\\
 & 500 &  4.261 & 142.6  & 1162.3\\
\hline
 \multirow{4}{*}{20} & 10 & 0.405  &  103.8  & 2015.4\\
 & 20 &  1.121 & 167.5  & 2906.7\\
 & 50 &  88.371 & 1568.9  & 127030.7 \\
 & 100 &  154.587 & 2257.2  & 122514.8 \\
 & 500 & 277.127  & 1800.5  & 111079.2\\
 \hline
\end{tabular}
\label{res:entropy_min}
\end{center}
\end{table}

To this point, we have demonstrated the modeling power of the multi-set extension of submodularity---$k$-submodularity---in real-life applications. These decision-making problems give rise to nonlinear optimization problems over $k$-sets that cannot be directly entered into optimization solvers. We have discussed the ways to reformulate these challenging problems as mixed-integer linear programs. The DCG algorithm leverages the power of MILP solvers and efficiently solves the reformulated problems to optimality despite the large numbers of constraints. Next, we provide a brief overview of another notion of generalized submodularity.

\section{Mixed-Integer Extension of Submodularity} 
\label{sect:mi}
Recall in Section \ref{sect:introduction}, we have motivated the mixed-integer extension of submodularity, in which 
we are allowed to place multiple, and even fractional, copies of the provided homogenous items in each $\lq$bin' (see Figure \ref{fig:sub_extensions}). Instead of having binary decision variables as in submodular set function optimization, now we need to incorporate mixed-integer decision variables. Specifically, we focus on \emph{Diminishing Returns (DR)-submodularity} in this tutorial.

\begin{figure}
  \centering 
            \includegraphics[width=6cm]{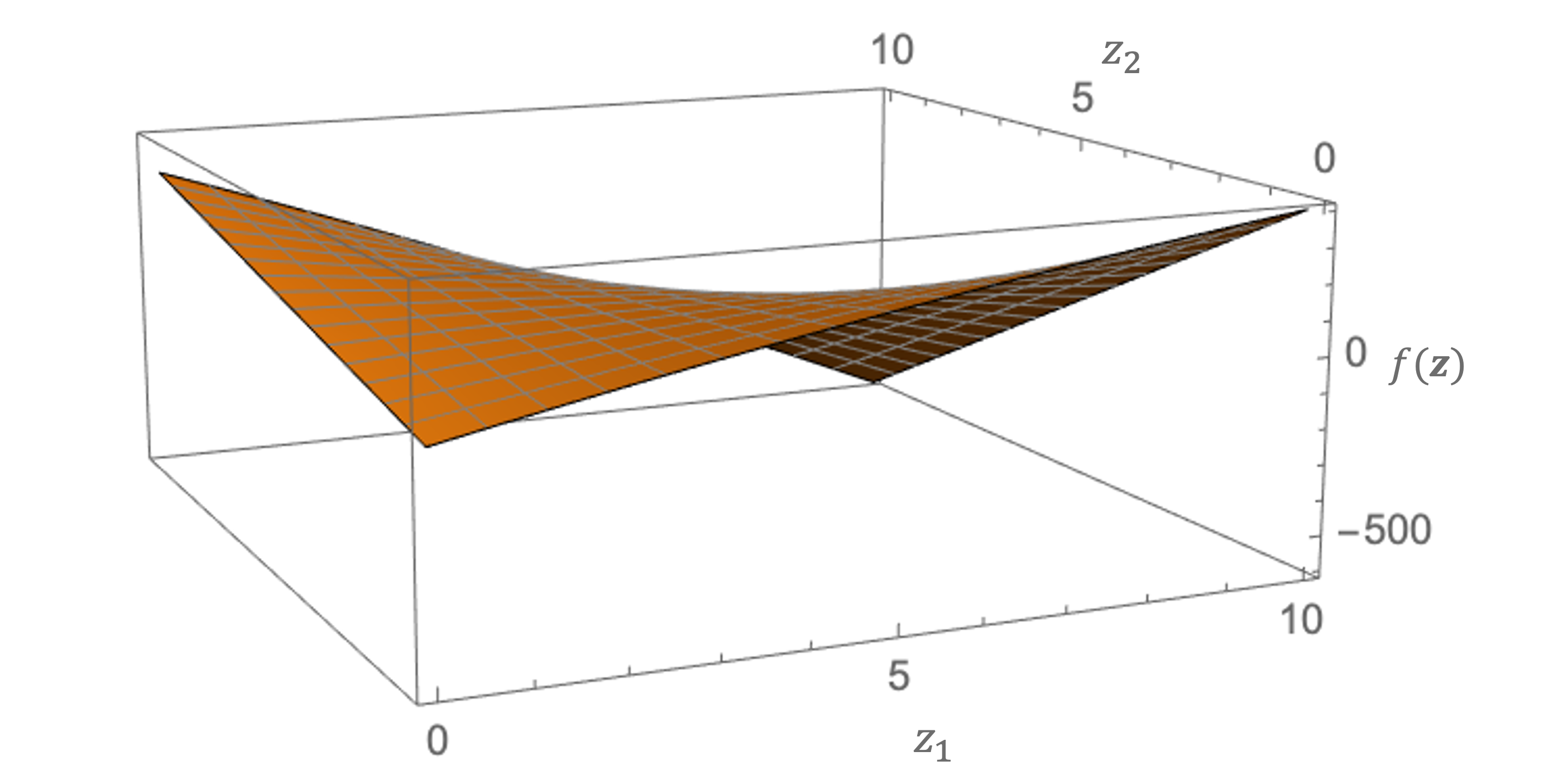}
            \includegraphics[width=6cm]{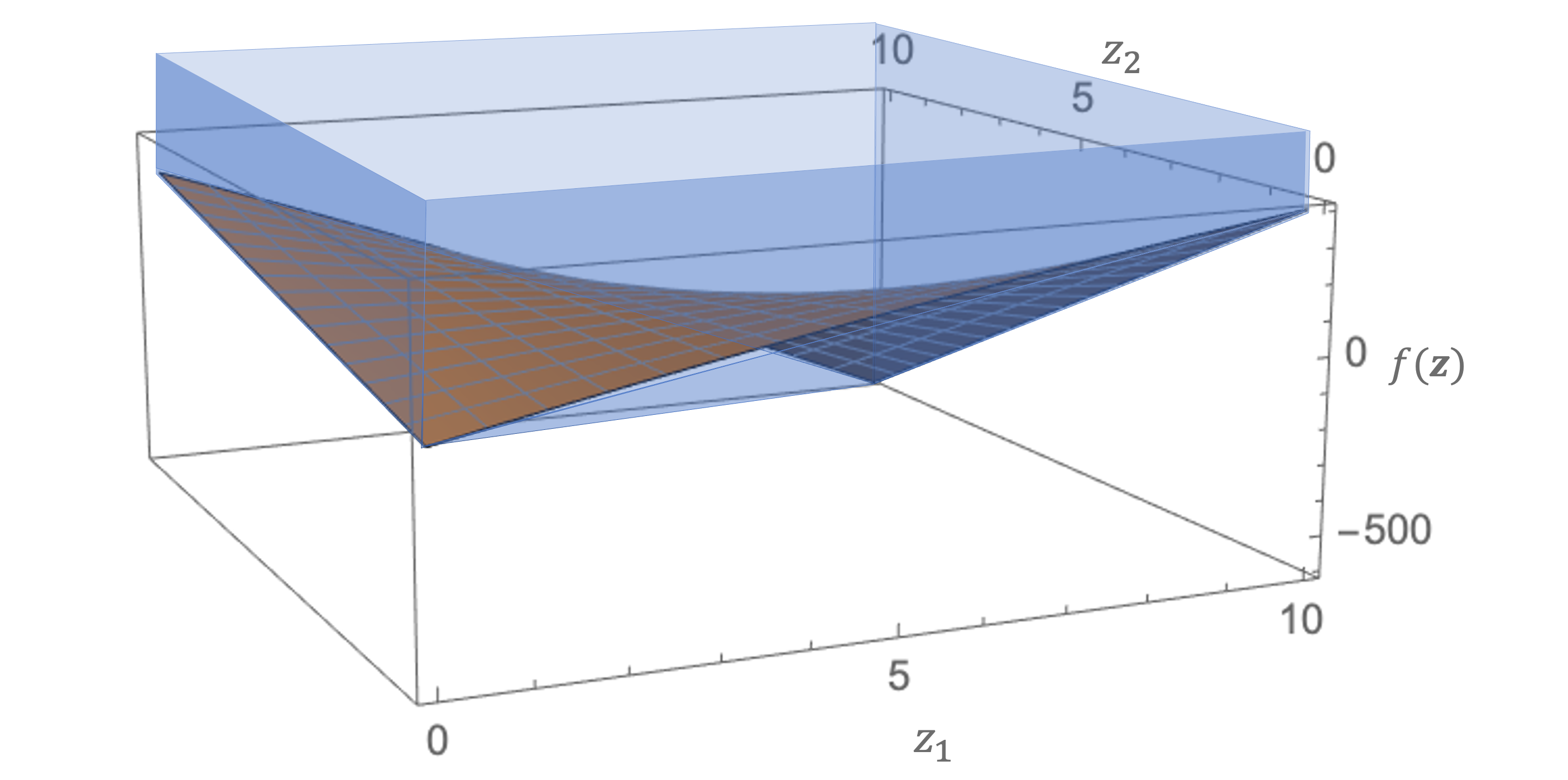}
           \caption{\small A continuous nonconvex and nonconcave DR-submodular function $f(\mathbf{x}) = -x_1^2 - 13x_1x_2+50x_1+30x_2$ and the convex hull of its epigraph under box constraints. The convex hull is polyhedral. }
           \label{fig:DR_nonconvex}
  \end{figure}

Let $\mathbf{e}^i\in\mathbb{R}^n$ be a vector with one in the $i$-th entry and zero everywhere else. 
\begin{definition}
\label{def:dr}
A function $f:\mathcal{X}\subseteq \mathbb{R}^n\rightarrow \mathbb{R}$ is \emph{DR-submodular} if 
\[f(\mathbf{x}+\alpha\mathbf{e}^i)-f(\mathbf{x})\geq f(\mathbf{y}+\alpha\mathbf{e}^i)-f(\mathbf{y})\] 
for every $i\in \{1,2,\dots, n\}$, for all $\mathbf{x},\mathbf{y}\in\mathcal{X}$ with $\mathbf{x}\leq\mathbf{y}$ component-wise, and for all $\alpha\in\mathbb{R}_+$ such that $\mathbf{x}+\alpha\mathbf{e}^i, \mathbf{y}+\alpha\mathbf{e}^i\in\mathcal{X}$. 
\end{definition}
Similar to Definition \ref{def:submodular_equiv} for submodular set functions, Definition \ref{def:dr} reflects the diminishing returns intuition. The marginal contribution of adding a positive step $\alpha$ in the $i$-th coordinate to a component-wise smaller vector $\mathbf{x}$ is higher than that to a larger vector, $\mathbf{y}$. A continuous DR-submodular function can be convex, concave, or neither. Figure \ref{fig:DR_nonconvex} (left) is an example of a nonconvex and nonconcave DR-submodular function. A twice differentiable function is DR-submodular if and only if its Hessian entries are nonpositive at every point in the domain \cite{bian2017guaranteed}. 

There exists another related mixed-integer extension of submodularity. 
\begin{definition}
\label{def:lattice}
Consider a function $f:\mathcal{X}\subseteq \mathbb{R}^n \rightarrow\mathbb{R}$. Suppose 
\[f(\mathbf{x}) + f(\mathbf{y}) \geq f(\mathbf{x}\wedge \mathbf{y}) + f(\mathbf{x}\vee \mathbf{y}),\]  
for all $\mathbf{x},\mathbf{y}\in\mathcal{X}$, where $\wedge$ and $\vee$ are the component-wise minimum and maximum operators, respectively. Then $f$ is referred to as \emph{submodular} (see \cite{bach2019submodular, staib2017robust, bian2017guaranteed, topkis1978minimizing, han2022polynomial}), or \emph{lattice submodular} when $\mathcal{X}\subseteq \mathbb{Z}^n$ (see \cite{ene2016reduction, soma2015generalization, soma2018maximizing}).
\end{definition}
Note that DR-submodularity is a subclass of (lattice) submodularity \cite{bian2017guaranteed}. Due to its diminishing returns property, DR-submodularity is prevalent in many applications, including but not limited to, optimal budget allocation, revenue management, sensor energy management, stability number of graphs, and mean field inference \cite{bian2017guaranteed, motzkin1965maxima, soma2017non, soma2018maximizing, sadeghi2021faster}. Therefore, we focus on DR-submodularity in this section.

\subsection{Applications of DR-submodularity}
We now provide detailed descriptions of two applications of DR-submodular functions. 

\subsubsection{Advertisement planning.} We are interested in designing an advertisement plan with optimal influence. Let $N = \{1,2,\dots, n\}$ be a set of media platforms (e.g., inline ads, TV commercial spots, social media influencers) and $T$ be the set of our target audience. Each media option $i\in N$ has a probability $p_{it}$ of influencing the audience $t\in T$. Our decisions are the ad volumes $x_i$ via media option $i$, for all $i\in N$. This variable could be continuous (e.g., TV air time) or discrete (e.g., number of social media influencers). Given a plan $\mathbf{x}\in\mathbb{R}^n$, the total influence of this plan is $f(\mathbf{x}) = \sum_{t\in T} \left[1 - \prod_{i\in N} (1-p_{it})^{x_i}\right]$, which is DR-submodular \cite{bian2017guaranteed}.

\subsubsection{Sensor energy management.} \label{eg:sensor_energy} The global energy crisis has made energy efficiency a crucial consideration in many applications. We describe the problem of placing sensors with adjustable energy levels to detect contamination events \cite{soma2015generalization}. We are given $N = \{1,2,\dots, n\}$ candidate locations for sensor placement, and our task is to decide the level of energy $x_i$ made available to $i\in N$. The amount of energy could fall under a continuous range or discrete levels, making our decision variables $\mathbf{x}$ continuous and discrete, respectively. We let $x_i=0$ if no sensor is placed in $i$.  When a contamination event $e\in E$ occurs, $z(i, e)$ represents the time taken for the pollution to reach sensor location $i$ in this event. Let $0< p < 1$ be a constant representing the chance of detection failure per unit of energy at any sensor location. When the pollution reaches sensor location $i$, the probability of successfully detecting it is $1-p^{x_i}$. We note that higher energy increases the chance of successful detections. Given $e\in E$, let $I_e$ be a random variable representing the first sensor that detects the pollution in event $e$. Let $z_\infty = \max_{i\in N, e\in E} z(i, e)$ be the maximal time for the pollution in event $e$ to reach any sensor location $i\in N$. If no sensor detects the pollution in event $e$, then we let $z(i, e) = z_\infty$. The expected time reduction in pollution detection is modeled by the function $f(\mathbf{x}) = \mathbb{E}_{e\in E} \mathbb{E}_{I_e} [z_\infty - z(I_e, e)]$. It has been shown that function $f$ is DR-submodular \cite{soma2015generalization}.

\subsection{DR-submodular optimization}
Motivated by numerous applications, DR-submodular \emph{maximization} has attracted much attention. For \emph{continuous} DR-submodular maximization, researchers have developed approximation and global optimization algorithms \cite{bian2017guaranteed, bian2017continuous, hassani2017gradient, sadeghi2021faster, medal2022spatial}. Other studies examine DR-submodular maximization with \emph{integer decision variables} and propose approximation algorithms with theoretical guarantees \cite{soma2017non, soma2018maximizing}. 

Few studies have focused on DR-submodular \emph{minimization}. They generally assume the decision variables to be box-constrained (i.e., $0 \leq x_i \leq u_i$ for $u_i\in\mathbb{R}_+$, $i\in N$). In \cite{ene2016reduction}, the authors reduce DR-submodular optimization with \emph{integer decision variables} to submodular set function optimization, which implies that integer DR-submodular minimization under box-constraints is polynomially solvable. The complexity depends on the size of the decision space and the logarithm of the  parameter values of the box-constraints (see also \cite{staib2017robust}). 

Recall that (lattice) submodularity subsumes DR-submodularity, so the works on (lattice) submodular minimization also apply to DR-submodular minimization. Topkis \cite{topkis1978minimizing} explores the properties of the minimizers for a class of parametric (lattice) submodular minimization problems. Bach \cite{bach2019submodular} generalizes Choquet integral (equivalent to Lov\`asz extension) to (lattice) submodular functions and proposes an algorithm for minimizing (lattice) submodular functions under box constraints. The time complexity of the algorithm is polynomial in the decision space size and the box constraints parameters and is thus pseudo-polynomial. Han et al. \cite{han2022polynomial} consider a class of constrained submodular minimization problems with \emph{mixed-binary} variables and establish its polynomial solvability. We take a polyhedral approach to tackle constrained DR-submodular minimization problems with \emph{mixed-integer} variables in \cite{yu2022drsubmod}.

\subsubsection{DR-submodular minimization}
We next discuss how the mixed-integer programming approach is applied to tackle constrained mixed-integer DR-submodular minimization \cite{yu2022drsubmod}. 

Consider the following problem: 
\begin{equation}
\label{eq:DR_min}
\min_{\mathbf{x}\in\mathcal{X}(\mathcal{G},\mathbf{u})} f(\mathbf{x}), 
\end{equation}
where $\mathcal{X}(\mathcal{G},\mathbf{u})$ is defined by box constraints and possibly monotonicity constraints: 
\begin{equation}
\label{eq:ZGU}
\mathcal{X}(\mathcal{G},\mathbf{u}) := \{\mathbf{x}\in\mathbb{Z}^n\times\mathbb{R}^m : \mathbf{0}\leq \mathbf{x}\leq \mathbf{u}, x_i \leq x_j, \: \forall\: (i,j)\in\mathcal{A}\}.
\end{equation}
Our decision $\mathbf{x}$ contains $n \in \mathbb Z_+$ general integer variables ($n\ge 1$) and $m \in \mathbb Z_+$ continuous variables. We let $\mathcal{V} = \{1,\dots, n+m\}$ be the index set of all the variables, and let $N$ be the set of indices corresponding to the integer variables. Notation $\mathcal{G}$ refers to a digraph with vertices $\mathcal{V}$ and arcs $\mathcal{A}$. The arc set $\mathcal{A}\subset \mathcal{V}\times \mathcal{V}$ governs the monotonicity constraints for the feasible set. The monotonicity constraints are inspired by monotone systems of linear inequalities, or the inequalities each with two variables and coefficients of opposite signs \cite{cohen1991improved, shostak1981deciding, aspvall1980polynomial, hochbaum1994simple}. They naturally arise in applications as well. For example, in sensor energy management (see Section \ref{eg:sensor_energy}), certain sensor locations provide more helpful information than others, thus we may be interested in imposing a partial ordering on the energy allocation. 

\begin{figure}
\centering
\includegraphics[width=5cm]{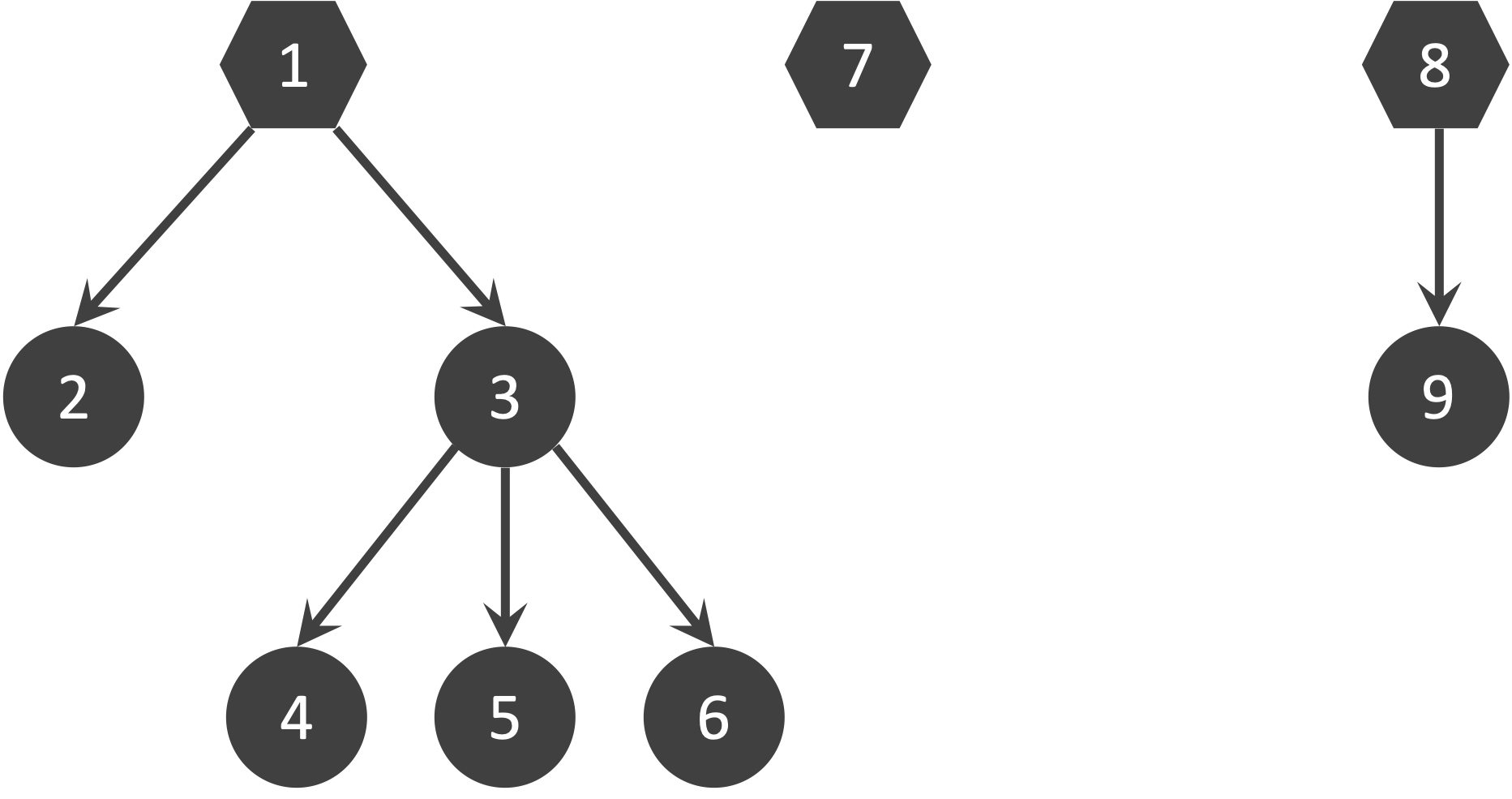}
\caption{\small An example of a directed rooted forest. }
\label{fig:DRF}
\end{figure}

We allow $\mathcal{G}$ to be any directed rooted forest, which is a disjoint union of directed trees, each with a designated root and all arcs pointing away from the root. Figure \ref{fig:DRF} is an example of a directed rooted forest. Such digraphs do not contain cycles, which would imply redundant equalities. We set forth the following definitions related to a directed rooted forest. For $i\in \mathcal{V}$, {we use  $R^+(i)\subseteq \mathcal{V}$ to denote the vertices that $i$ can reach following the arcs in $\mathcal{A}$. We call these vertices the \emph{descendants} of $i$. }

Surprisingly, the constrained mixed-integer nonlinear optimization problem \eqref{eq:DR_min} can be reformulated as a continuous linear program. To see this, we first rewrite problem  \eqref{eq:DR_min} as 
\begin{equation}
\label{eq:DR_min_epi}
\min \left\{w: (\mathbf{x},w)\in\conv{\mathcal{P}_f^{\mathcal{X}(\mathcal{G},\mathbf{u})}}\right\}, 
\end{equation}
where 
\[\mathcal{P}_f^{\mathcal{X}(\mathcal{G},\mathbf{u})} := \{(\mathbf{x},w)\in\mathcal{X}(\mathcal{G},\mathbf{u})\times\mathbb{R}: w\geq f(\mathbf{x})\}\]
is the epigraph of $f$ under $\mathcal{X}(\mathcal{G},\mathbf{u})$. Theorem \ref{thm:tree_m_conv} proposes a full characterization of $\conv{\mathcal{P}_f^{\mathcal{X}(\mathcal{G},\mathbf{u})}}$, and this convex hull happens to be polyhedral. Therefore, \eqref{eq:DR_min_epi} becomes a linear program with continuous variables and exponentially many linear constraints. Figure \ref{fig:DR_nonconvex} (right) pictorially illustrates this result on an example. 

\begin{theorem} \cite{yu2022drsubmod} 
\label{thm:tree_m_conv}
All the DR inequalities associated with valid permutations, along with the Mixed-Integer Rounding (MIR) inequalities, the box constraints, and the monotonicity constraints that define $\conv{\mathcal{X}(\mathcal{G},\mathbf{u})}$, fully describe $\conv{\mathcal{P}_f^{\mathcal{X}(\mathcal{G},\mathbf{u})}}$, under some technical conditions.  
\end{theorem}

We next define the MIR inequalities and the DR inequalities mentioned in Theorem \ref{thm:tree_m_conv}. The MIR inequalities are needed to convexify the mixed-integer feasible set $\mathcal{X}(\mathcal{G},\mathbf{u})$. Let 
\[\Psi := \{\psi\in \mathcal{V} : u_\psi \notin\mathbb{Z}, R^+(\psi)\cap N\neq \emptyset\}\] 
be  a collection of the vertices corresponding to the non-integer upper-bounded variables that each has at least one discrete descendant. For any $\psi\in\Psi$, $u_\psi\notin\mathbb{Z}$. The existence of such variables tends to make the continuous relaxation of $\mathcal{X}(\mathcal{G},\mathbf{u})$ differ from its convex hull. Thus, we need additional valid inequalities to strengthen the continuous relaxation. The MIR inequality associated with $\psi$ is the following linear inequality: 
 \[-x_{\ch{\psi}} - \frac{u_\psi - x_\psi}{u_\psi - \floor{u_\psi}} \leq -\ceil{u_\psi}, \]
 {where $\ch{\psi}$ is a vertex with $(\psi, \ch{\psi}) \in\mathcal{A}$. We provide an example of MIR inequalities in Example \ref{eg:MIR}. }
 \begin{example}
 \label{eg:MIR}
 Consider the feasible set $\mathcal{X}(\mathcal{G},\mathbf{u}) = \{x_1\in\mathbb{R}, x_2\in\mathbb{Z}: 0 \leq x_1 \leq 2.4, 0\leq x_2 \leq 3, x_1 \leq x_2\}$. We notice that $\Psi  = \{1\}$. In Figure \ref{fig:mir}, the solid black dot and lines represent the mixed-integer feasible points. The dotted region is the continuous relaxation of $\mathcal{X}(\mathcal{G},\mathbf{u})$, and the purple region highlights $\conv{\mathcal{X}(\mathcal{G},\mathbf{u})}$. We observe that the continuous relaxation of $\mathcal{X}(\mathcal{G},\mathbf{u})$ does not match the convex hull. The dashed purple line corresponds to the MIR inequality $-x_2 - (2.4 - x_1)/0.4 \leq -3$ associated with $\psi = 1$. This MIR inequality and the box and monotonicity constraints fully characterize $\conv{\mathcal{X}(\mathcal{G},\mathbf{u})}$.
\begin{figure}[h] 
   \centering
    \includegraphics[width=11cm]{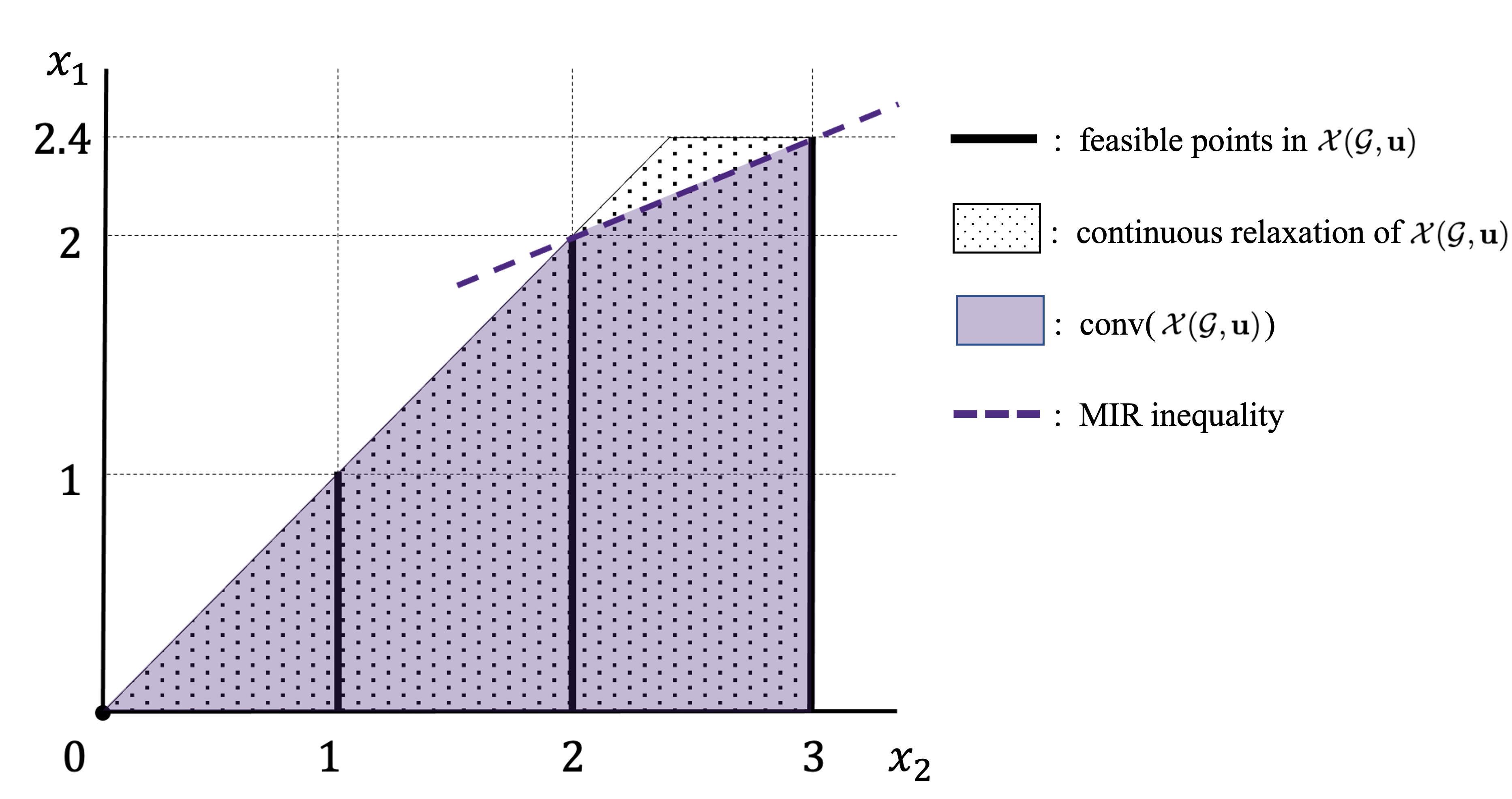}
    \caption{\small An example of an MIR inequality. }
    \label{fig:mir}
\end{figure}
 \end{example}
As stated in the next theorem, the MIR inequalities associated with all $\psi\in\Psi$ are important for the convex hull description of the mixed-integer feasible set. 
 
 \begin{theorem} \cite{yu2022drsubmod}
The MIR inequalities associated with all $\psi\in\Psi$, and the box and monotonicity constraints, sufficiently describe the convex hull of the feasible set $\mathcal{X}(\mathcal{G},\mathbf{u})$, under some technical conditions. 
 \end{theorem}
 
The DR inequalities further involve the epigraph variable $w$. Let $\bfs{\delta}$ be a permutation of $\mathcal{V}$. We define a \emph{DR inequality} associated with $\bfs{\delta}$ by
\begin{equation}
\label{eq:DR}
w \geq \sum_{k=1}^{|\mathcal{V}|} t(\bfs{\delta},\mathbf{x})_k[f(P(\bfs{\delta},k)) - f(P(\bfs{\delta},k-1))],
\end{equation}

For every $k\in \{0,1,\dots,|\mathcal{V}|\}$, $t(\bfs{\delta},\mathbf{x})_k$ is a linear function of $\mathbf{x}$, and $P(\bfs{\delta},k)$ is an extreme point of $\conv{\mathcal{X}(\mathcal{G},\mathbf{u})}$.  The closed forms of $t(\cdot, \cdot)$ and $P(\cdot)$ are provided in \cite{yu2022drsubmod}. We remark that all DR inequalities are linear and homogenous. Furthermore, they reduce to the celebrated EPIs \eqref{eq:EPI} when $f$ is a classical submodular function. Given that every DR inequality is associated a permutation of $\mathcal{V}$ (that satisfies certain properties), there are exponentially many DR inequalities. Each DR inequality can be separated exactly using an $\mathcal{O}(|\mathcal{V}|^2\log |\mathcal{V}|)$ algorithm  \cite{yu2022drsubmod}. Due to the equivalence between separation and optimization, we reach the following conclusion: 

\begin{corollary} \cite{yu2022drsubmod}
DR-submodular minimization over the mixed-integer feasible set $\mathcal{X}(\mathcal{G},\mathbf{u})$ is polynomial-time solvable, under some technical conditions. 
\end{corollary}

\section{Conclusions}
\label{sect:conclusion}
In this tutorial, we introduce generalized submodularity by providing an overview on submodularity and two of its extensions. We show the modeling power of generalized submodularity by describing examples of its applications in detail. These applications give rise to mixed-integer nonlinear optimization problems with generalized submodular objectives, and these optimization problems are challenging to solve exactly, due to the discrete decision space and the nonlinearity. We overcome these challenges with mixed-integer programming techniques. Specifically, we derive mixed 0-1 or continuous linear programming reformulations of the original problems and devise effective exact solution methods by applying the Delayed Constraint Generation framework and leveraging the power of mixed-integer linear programming solvers.

\section*{Acknowledgement}
The research is supported, in part, by ONR Grant N00014-22-1-2602. 

\bibliography{gso_ref}{}
\bibliographystyle{apalike}

\end{document}